\documentclass[a4paper,12pt, reqno]{amsart}

\usepackage[svgnames]{xcolor}
\usepackage[verbose]{geometry}
\usepackage[bbgreekl]{mathbbol}

\usepackage{amsmath, amssymb, amsthm, enumerate, url}
\usepackage{euscript}
\usepackage{amsfonts}
\usepackage{enumitem}
\usepackage{comment}
\usepackage{tikz-cd}
\usepackage{quiver}
\usepackage{mathrsfs}
\usepackage{mathtools}
\usepackage{graphicx}
\usepackage{thmtools}
\usepackage{microtype}
\usepackage{braket}
\renewcommand{\Set}{\set}
\let\Set\relax

\usepackage{stmaryrd}
\SetSymbolFont{stmry}{bold}{U}{stmry}{m}{n}

\usepackage[backend=biber, style=trad-alpha, doi=true, isbn=false, url=true, date=year]{biblatex}
\usepackage[svgnames]{xcolor}
\definecolor{darkgreen}{rgb}{0,0.45,0}
\definecolor{darkred}{rgb}{0.75,0,0}
\definecolor{darkblue}{rgb}{0,0,0.6}
\definecolor{acid}{HTML}{438806}
\definecolor{crimson}{RGB}{177,12,12}
\definecolor{teal}{HTML}{007B7F}
\definecolor{orange}{HTML}{F28C28}
\definecolor{coral}{HTML}{FF6F61}
\definecolor{plum}{HTML}{804568}
\definecolor{blue(munsell)}{rgb}{0.0, 0.5, 0.69}
\usepackage{hyperref}
\hypersetup{
	colorlinks = true,
	urlcolor=plum,
	linkcolor=teal,
	citecolor=acid,
}

\usepackage{cleveref}
\Crefname{ex}{Example}{Examples}
\Crefname{defn}{Definition}{Definitions}

\DeclareSourcemap{
  \maps[datatype=bibtex]{
    \map[overwrite]{
      \step[fieldsource=url,
        match=\regexp{\Ahttps?://(?:www\.)?arxiv\.org/(?:abs|pdf)/},
        final]
      \step[fieldsource=url, fieldtarget=eprint]
      \step[fieldsource=eprint,
        match=\regexp{\Ahttps?://(?:www\.)?arxiv\.org/(?:abs|pdf)/},
        replace={}]
      \step[fieldsource=eprint,
        match=\regexp{\.pdf\z}, replace={}]
      \step[fieldset=eprinttype, fieldvalue={arxiv}]
    }
  }
}

\DeclareFieldFormat{eprint:arxiv}{%
  Available at\addspace
  \href{https://arxiv.org/abs/#1}{arXiv\addcolon#1}}

\DeclareFieldFormat{url}{Available at\addspace\url{#1}}

\usepackage{soul}
\setul{3pt}{.4pt}
\usepackage{dsfont}

\usetikzlibrary{arrows.meta}

\usepackage{tikz}
\usepackage{tikz-cd}

\usetikzlibrary{calc}
\usetikzlibrary{fit}
\usetikzlibrary{shapes}
\usetikzlibrary{arrows}
\usetikzlibrary{decorations.markings}
\usetikzlibrary{decorations.pathmorphing}
\usetikzlibrary{positioning}
\usetikzlibrary{intersections}

\tikzset{
  commutative diagrams/arrow style=tikz,
}

\tikzset{cd-style/.style={commutative diagrams/every diagram}}
\tikzset{cd-arrow-style/.style={commutative diagrams/.cd, every arrow, every label}}

\makeatletter
\tikzset{
  label/.style n args={2}{
    edge node={node [
      execute at begin node=\iftikzcd@mathmode$\fi,
      execute at end node=\iftikzcd@mathmode$\fi,
      /tikz/commutative diagrams/.cd,every label,
      #2
      ] {#1}}
  }
}
\makeatother

\tikzset{fibtip/.tip={Triangle[open,angle=60:4.5pt]}}
\tikzset{tfibtip/.tip={Bar[sep]Triangle[open,angle=45:4pt]}} 

\pgfarrowsdeclaredouble{fibbtip}{fibbtip}{fibtip}{.fibtip}

\tikzset{ulttip/.tip={Computer Modern Rightarrow[reversed]}}

\pgfarrowsdeclarealias{c}{c}{left hook}{right hook}
\pgfarrowsdeclarealias{c'}{c'}{right hook}{left hook}

\tikzset{fib/.code={\pgfsetarrowsend{fibtip}}}
\tikzset{fibb/.code={\pgfsetarrowsend{fibbtip}}}
\tikzset{tfib/.code={\pgfsetarrowsend{tfibtip}}}
\tikzset{inj/.code={\pgfsetarrowsstart{c}}}
\tikzset{inj'/.code={\pgfsetarrowsstart{c'}}}
\tikzset{cover/.style={->>}}
\tikzset{tcof/.style={tail}}
\tikzset{zigzag/.style={commutative diagrams/rightsquigarrow}}
\tikzset{ulttip/.code={\pgfsetarrowsend{ulttip}}}

\tikzset{suplabel/.style={label={#1}{auto=left,pos=1}}}
\tikzset{sublabel/.style={label={#1}{auto=right,pos=1}}}
\tikzset{sub/.style={label={#1}{auto=right,pos=1}}}

\tikzset{lw/.style={"lw"{#1}},lw/.default={very near end}}
\tikzset{lw'/.style={"lw"'{#1}},lw'/.default={very near end}}
\tikzset{weq/.style={"\sim"}}
\tikzset{weq'/.style={"\sim"'}}
\tikzset{ult/.style={ulttip,"#1"'{very near end}}}
\tikzset{ult'/.style={ulttip,"#1"{very near end}}}

\newcommand{\generalto}[2]{ \mathrel{\mkern-1mu
  \tikz[baseline={([yshift=-0.58ex]a.south)}]{%
    \node[minimum width=1em,align=center,inner xsep=0.5ex,inner ysep=0.25ex] (a) {$\scriptstyle #2$};
    \draw[cd-arrow-style,#1] (a.south west) -- (a.south east);}
 \mkern-1mu}}

\renewcommand{\to}[1][]{ \generalto{->}{#1} }

\newcommand{\coverto}[1][]{ \generalto{->>}{#1\ \,} }
\newcommand{\vuto}[2][]{ \mathrel{{\generalto{ulttip}{#1}}\mkern2mu_{{#2}}\mkern-2mu}}
\newcommand{\vuleq}[2][]{ \mathrel{\preccurlyeq_{{#2}}^{{#1}}\mkern-2mu}}

\setlist[enumerate]{label={(\roman*)}, ref={(\roman*)}} 

\theoremstyle{plain}
\newtheorem{Thm}{Theorem}[subsection]
\newtheorem*{nThm}{Theorem}
\newtheorem{prop}[Thm]{Proposition}
\newtheorem{cor}[Thm]{Corollary}
\newtheorem{lem}[Thm]{Lemma}
\theoremstyle{definition}
\newtheorem{defn}[Thm]{Definition}

\newtheorem{cons}[Thm]{Construction}
\newtheorem{nota}[Thm]{Notation}
\newtheorem*{nnota}{Notation}
\newtheorem{ex}[Thm]{Example}
\newtheorem{rmk}[Thm]{Remark}

\DeclareMathOperator{\dom}{dom}

\DeclareMathOperator*{\colim}{colim}
\DeclareMathOperator{\cl}{\mathsf{cl}}

\newcommand{\mono}{\hookrightarrow}
\newcommand{\epi}{\twoheadrightarrow}

\def\slashedrightarrow{\relbar\mathrel{\mkern-4mu}\joinrel\mapstochar\mathrel{\mkern-4mu}\joinrel\rightarrow}
\newcommand{\pro}{\slashedrightarrow}

\newcommand{\E}{\mathscr{E}}

\newcommand{\X}{\mathscr{X}}

\newcommand{\Z}{\mathscr{Z}}
\newcommand{\Y}{\mathscr{Y}}

\newcommand{\ptuf}{1}

\newcommand{\op}{\mathrm{op}}

\newcommand{\CAT}{\mathsf{CAT}}
\newcommand{\TopSp}{\mathsf{TopSp}}

\newcommand{\pt}{\mathsf{pt}}

\newcommand{\id}{\mathsf{i}}
\newcommand{\diag}{\mathsf{d}}

\newcommand{\Set}{\mathsf{Set}}

\newcommand{\UF}{\mathds{U}}

\newcommand{\vuCAT}{\mathsf{vuCAT}}
\newcommand{\vuPos}{\mathsf{vuPos}}

\newcommand{\mrg}{\curlyvee}
\newcommand{\bigmrg}{\bigcurlyvee}

\newcommand{\Tref}[1]{|#1|}

\newcommand{\FF}{\mathds{F}}

\newcommand{\vfCAT}{\mathsf{vfCAT}}
\newcommand{\vfPos}{\mathsf{vfPos}}

\DeclareFontFamily{U}{dmjhira}{}
\DeclareFontShape{U}{dmjhira}{m}{n}{ <-> dmjhira }{}

\makeatletter
\let\old@tocwrite\@tocwrite
\def\@tocwrite#1#2{\ifx\@secnumber\@empty\else\old@tocwrite{#1}{#2}\fi}
\makeatother

\newcommand{\emaillink}[1]{\href{mailto:#1}{\sf #1}}

\title[Categorified convergence]{%
  Categorified convergence:\\
  \makebox[\linewidth][c]{%
    \small dropping the `ultra' from `virtual ultracategories'%
  }%
}

\date{\today}

\author{Gabriel Saadia}
\author{Umberto Tarantino}

\address{
Gabriel \textsc{Saadia} \newline
Department of Mathematics\newline
Stockholm University\newline
Stockholm, Sweden\newline
\emaillink{gabriel.saadia@math.su.se}
}

\address{
Umberto \textsc{Tarantino} \newline
Université Paris Cité, CNRS, IRIF, F-75013\newline
Paris, France\newline
\href{mailto:tarantino@irif.fr}{\sf tarantino@irif.fr}
}

\thanks{The first-named author was supported by the Knut and Alice Wallenberg Foundation (KAW) Foundation, under the project ``Type Theory for Mathematics and Computer Science'' (PI Thierry Coquand). The second-named author acknowledges financial support from the Agence Nationale de la Recherche
(ANR), project ANR-23-CE48-0012-01, and from the Laboratoire Méthodes Formelles (LMF), Université Paris-Saclay.}

\begin{document}

\begin{abstract}
\looseness=-1

Virtual ultracategories are a recent categorification of topological spaces seen as sets of points together with a convergence relation between points and ultrafilters. However, a topology can also be described in terms of convergence of filters: in this paper, we thus introduce an analogue of virtual ultracategories based on filters rather than ultrafilters. At the topological level, convergence of filters is completely determined by its restriction to ultrafilters. Our main result shows that, under the assumption of tautness, the same phenomenon holds at the categorified level, thereby offering a non-\emph{ultra-biased} perspective on virtual ultracategories.

\smallskip \noindent \textbf{Keywords.} convergence, virtual ultracategory, ultrafilter, filter, ultraproduct, reduced product, topological space

\noindent  \textbf{MSC2020.} 54A20, 03C20, 18N25, 18F60, 18F20, 54B30, 03G30.
\end{abstract}

\maketitle

\vspace{-1em}

\setcounter{tocdepth}{2} 
\tableofcontents

\section*{Introduction}\label{sec:intro}
\addcontentsline{toc}{section}{\nameref{sec:intro}}

\subsection*{Topology via convergence}Topology has long been studied in terms of convergence. While this point of view was already implicit in the 19th-century study of $\mathbb N$-indexed sequences and their limits in the context of real analysis, a first systematic approach to topology via convergence is due to Moore and Smith \cite{mooreGeneralTheoryLimits1922}, who considered generalized sequences indexed by arbitrary directed sets. Their work was subsequently developed into the modern theory of \emph{nets} by Birkhoff \cite{birkhoffMooreSmithConvergenceGeneral1937}, Tukey \cite{tukeyConvergenceUniformityTopology1940}, and Kelley \cite{kelleyConvergenceTopology1950,kelleyGeneralTopology1955}. 
Convergence of nets can be equivalently formulated, using the language of \emph{filters} introduced by Cartan \cite{cartanTheorieFiltres1937,cartanFiltresUltrafiltres1937}, simply as a relation between filters on a set and points of the set. In these terms, a topology on a set $X$ is entirely characterized by its induced convergence relation: a filter on $X$ converges to a point $p$ precisely when it contains every neighborhood of $p$, and conversely a subset $U \subseteq X$ is open precisely when it belongs to every filter converging to a point in $U$. This observation initiated the study of convergence as a primitive notion, rather than merely a relation derived from a pre-existing topology, and it led to the development of abstract notions of `convergence spaces' generalizing the topological ones, notably in the works of Choquet \cite{choquetConvergence1948}, Kowalsky \cite{kowalskyLimesraumeUndKomplettierung1954}, Fischer \cite{fischerLimesraume1959}, and Kent \cite{kentConvergenceFunctionsTheir1964}.

\looseness=-1
The convergence relation of a topological space is, however, highly \emph{redundant}. Although already identified by Kelley in the language of nets (cf.\ the \emph{universal nets} of \cite{kelleyConvergenceTopology1950}), this redundancy becomes particularly evident in terms of filters: if a filter $\sigma$ converges to a point $p$, then so does every filter refining $\sigma$. In fact, classically, convergence is completely determined by its restriction to \emph{ultrafilters}, the maximal proper filters with respect to refinement:
\[ \sigma \mbox{ converges to }p \iff \mbox{every ultrafilter refining } \sigma \mbox{ converges to }p .\]
Thus, from the perspective of convergence, passing from filters to ultrafilters simply amounts to a change of presentation which loses no information about the topology of a space. The goal of this paper is to prove an analogue of this phenomenon for a suitable categorification of topological convergence relations appearing in categorical logic. 

\subsection*{Virtual ultracategories}\looseness=-1 In the three independent works \cite{saadiaExtendingConceptualCompleteness2026,hamadGeneralisedUltracategoriesConceptual2025,vangoolToposesEnoughPoints2026}, the ultrafilter-based presentation of topological spaces was categorified into the notion of \emph{virtual ultracategories}, henceforth called \emph{vu-categories}. Intuitively, a vu-category is a categorical structure where `objects' can be thought of as points and `arrows' can be thought of as witnessing convergence of indexed families of points with respect to ultrafilters on the index set.

As the name suggests, vu-categories also generalize Makkai's \emph{ultracategories} \cite{makkaiStoneDualityFirst1987}. The key idea underlying ultracategories is that, by \L o\'s's theorem, the class of models $\operatorname{Mod}(\mathbb T)$ of a finitary first-order theory $\mathbb T$ is closed under the formation of ultraproducts. When $\operatorname{Mod}(\mathbb T)$ is seen as an ultracategory, the ultraproduct of a family of models with respect to some ultrafilter then acts as a canonical `limit' of the family along that ultrafilter, much in the same way as ultrafilters admit unique limits in a compact Hausdorff space (cf.\ \cite{manesTripleTheoreticConstruction1969}). 

\looseness=-1
The key insight provided by vu-categories is that the existence of these canonical limits is not necessary to define a convergence structure on $\operatorname{Mod}(\mathbb T)$: if the latter is not closed under ultraproducts, we can still work with \emph{virtual} ultraproducts. This idea allowed \cite{saadiaExtendingConceptualCompleteness2026,hamadGeneralisedUltracategoriesConceptual2025,vangoolToposesEnoughPoints2026} to lift the reconstruction result of \cite{makkaiStoneDualityFirst1987,lurieUltracategories2018} from coherent logic to the setting of \emph{geometric} logic, where a \L o\'s-like theorem is not available.
However, once one is willing to work `virtually', the original motivation for restricting to ultraproducts disappears. 

\subsection*{Dropping the `ultra'}
The first observation at the origin of this work is that nothing in the definition of a vu-category inherently depends on ultrafilters, so that it can be formulated just as well in terms of arbitrary filters (cf.\ \cite[\S 5.4]{aristoteProfunctorialAlgebras2026}, and see also \cite[\S 6]{tarantinoUltracategoriesKanExtensions2025} for an analogous remark in the setting of ultracategories). Thus, our first contribution in this paper is to formally introduce \emph{vf-categories}, a variation of vu-categories categorifying the filter-based presentation of topological convergence relations. 

\looseness=-1
An important example of a vf-category is the class of models $\operatorname{Mod}(\mathbb T)$ of a geometric theory $\mathbb T$, whose convergence structure is described in complete analogy with the vu-categorical case. In this case, the role of ultraproducts is played by the more general \emph{reduced products}, well known in model theory (see, e.g., \cite[\S 6]{changModelTheory1990}). We emphasize here that such a shift would not be possible in the setting of ultracategories: the class of models of a first-order theory is generally not closed under reduced products. Instead, this issue vanishes at the `virtual' level, just as the lack of a \L o\'s-like theorem for geometric logic does in the context of vu-categories.

\looseness=-1 
Clearly, every vf-category yields a vu-category simply by forgetting the additional information provided by non-maximal filters: in this paper, we study the extent to which a vf-category is determined by its vu-categorical `trace'.
Our investigations will start in the topological setting, where we know that convergence of filters is determined by its restriction to ultrafilters. As remarked in \cite[Rem.\ 3.14]{vangoolToposesEnoughPoints2026} and \cite[Thm.\ 5.5]{hamadGeneralisedUltracategoriesConceptual2025}, topological spaces (seen as vu-categories) are characterized by an additional condition on the convergence relation, formulated and studied for arbitrary vu-categories in \cite{saadiaGeometricMorphismsVirtual2026} under the name of \emph{tautness}. Accordingly, we shall identify an appropriate tautness condition for vf-categories, allowing us to obtain a characterization of topological spaces among vf-categories. 
We then move to the general case of vf-categories and vu-categories: the main result of this paper is that, under the assumption of tautness, vf-categories and vu-categories describe the same notion of convergence. This result crystallizes taut vf-categories as a robust categorified analogue of topological convergence relations which is no longer \emph{ultra-biased}, paving the way towards a constructive exploration of the aforementioned reconstruction theorem.

\subsection*{Structure of the paper and contributions}\looseness=-1
We begin, in \Cref{sec:vf-cats}, by introducing the main characters of this story. In \S\ref{ssec:filters-ultrafilters} we recall the necessary background on filters, and in \S\ref{ssec:category-of-filters} we introduce their category. In \S\ref{ssec:completions}, we extend the construction of the category of filters to a pseudomonad acting on categories, which we employ in \S\ref{ssec:vf-cats-vu-cats} to define vf-categories (\Cref{defn:vf-cats}). Moreover, restricting to ultrafilters, we recall the definition of vu-categories from \cite{saadiaExtendingConceptualCompleteness2026,hamadGeneralisedUltracategoriesConceptual2025,vangoolToposesEnoughPoints2026}. In \S\ref{ssec:merging}, we define the \emph{merging} operation in a vf-category, and we conclude by proving some properties that we will use in the following sections.

In \Cref{sec:tautness}, we study the concept of \emph{tautness}. First, in \S\ref{ssec:tautness-definition}, we recall tautness for vu-categories from \cite{saadiaGeometricMorphismsVirtual2026} and we show how to adapt it to vf-categories. In \S\ref{ssec:tautness-sheaf}, we prove that tautness for vf-categories corresponds to a suitable sheaf condition (\Cref{cor:taut-neighbourhood}). As a consequence, in \S\ref{ssec:tautness-preserved} we show that tautness is preserved when passing from a vf-category to its underlying vu-category (\Cref{cor:U-preserves-tautness}); moreover, in \S\ref{ssec:taut-vf-posets} we characterize topological spaces as precisely the taut vf-posets (\Cref{cor:vf-poset=topsps}), which constitutes the posetal case of our main result. Finally, in \S\Cref{ssec:orthogonality}, we characterize tautness in terms of a useful orthogonality property (\Cref{prop:tautness-orthogonality-converse}).

In \Cref{sec:equivalence}, we prove the main theorem of the paper, stated below.

\begin{nThm}[\ref{thm:main-theorem}]
    The forgetful 2-functor $\vfCAT \to \vuCAT$ restricts to a 2-equivalence between taut vf-categories and taut vu-categories:
    \[\vfCAT_{\mathrm{taut}}\simeq\vuCAT_{\mathrm{taut}}.\]
\end{nThm}
\hspace{-\parindent}Concretely, our strategy will be to construct a taut vf-category out of a taut vu-category, and to show that this construction realizes an inverse to the forgetful 2-functor.\footnote{
We note here that our construction will require the starting vu-category to be taut: in other words, we will not obtain the above equivalence, say, by restricting a general adjunction between vf-categories and vu-categories. } 
To this end, the key technical ingredient is the notion of \emph{walking vf-arrows}, which we introduce in \S\ref{ssec:walking-vf-arrows}. Exploiting walking vf-arrows, we can endow a taut vu-category with the structure of a vf-category, which we assemble and which we prove to be taut in \S\ref{ssec:assembling} (\Cref{thm:obtaining-vf-cat} and \Cref{prop:FX-is-taut}). Finally, in \S\ref{ssec:main-theorem}, we prove the main theorem: we first show that this process allows us to reconstruct a taut vf-category from its underlying vu-category (\Cref{prop:hom-UF=id}), and then we conclude the proof.

\begin{nnota}
\looseness=-1 
    We fix here some notation and conventions used throughout the paper. First, all 2-categorical terminology and concepts employed in the paper are \emph{strict}: in particular, in a 2-category we speak of \emph{coproducts}, \emph{coequalizers}, and \emph{pushouts} to mean the strict 2-categorical notions of the same name. We then denote by $\CAT$ the usual 2-category of categories, always assumed to be \emph{locally small}. We denote by $\TopSp$ the 2-category of topological spaces and continuous maps with their pointwise specialization order.  We denote filters by Greek letters $\sigma, \tau, \lambda, \kappa$, their carrier sets by uppercase letters $S, T, L, K$, and their generic element by the corresponding lowercase letters $s,t, \ell, k$. For a set $S$ and an $S$-indexed family of sets $(X_s)_{s\in S}$, we denote by $S \cdot X_s$ the disjoint union $\coprod_{s\in S} X_s$: this perhaps unusual notation will be justified by \Cref{prop:F-is-lextensive}. Finally, the word \emph{topos} invariably means \emph{Grothendieck topos}. 
\end{nnota}

\subsection*{AI declaration}In this paper, AI was only used as an enhanced Google Search to retrieve existing literature, and to correct typos and grammatical mistakes. None of the ideas, statements, or results were obtained with AI assistance.

\section{Vf-categories}\label{sec:vf-cats}

\subsection{Filter combinatorics}\label{ssec:filters-ultrafilters}
We begin by recalling the necessary background on \emph{filters}, so as to fix our notation and conventions. 

\begin{defn}[Filters]
    \looseness=-1
    A filter $\sigma$ on a set $S$ is a nonempty set of subsets of $S$ which is upward-closed and closed under finite intersections; subsets of $S$ lying in $\sigma$ are said to be \emph{$\sigma$-large}. We say that a filter $\tau$ \emph{refines} $\sigma$ if $\sigma \subseteq \tau$, that is, if any $\sigma$-large set is $\tau$-large. A filter is \emph{proper} if it does not contain $\emptyset$. A proper filter is \emph{principal} if it contains a singleton.
\end{defn}

\begin{defn}[Ultrafilters]
    A proper filter is an \emph{ultrafilter} if it is maximal among proper filters with respect to the order given by refinement. Equivalently, a filter $\sigma$ on a set $S$ is an ultrafilter if for each subset $A\subseteq S$ exactly one of $A$ and $S\setminus A$ is $\sigma$-large. Principal filters are clearly ultrafilters.
\end{defn}

\begin{rmk}
    By the \emph{ultrafilter principle}, a weak form of the axiom of choice, every proper filter can be refined to an ultrafilter. The \emph{Fréchet filter} on a set is the filter given by its cofinite subsets; on an infinite set, every nonprincipal ultrafilter refines the Fréchet filter.
\end{rmk}

\begin{nota}
    For a filter $\sigma$ on a set $S$, if $P$ is a property depending on $s \in S$, we write 
    \[ P(s) \mbox{ for } s:\sigma \quad \mbox{ or, in symbols, }\quad \forall_{s:\sigma} P(s) \]
    if $P(s)$ holds for all $s$ in a $\sigma$-large subset of $S$. 
\end{nota}

\begin{defn}[Filter and ultrafilter monads]\label{def:filter-monads}
    Mapping each set $S$ to the set $\gamma S$ of filters on it determines the \emph{filter monad} $\gamma \colon \Set \to \Set$, where for a function $\phi \colon T \to S$ and a filter $\tau$ on $T$:
    \[\gamma (\phi)(\tau) \coloneqq \set{ S_0 \subseteq S | \phi^{-1}(S_0) \in \tau}, \text{ also denoted by } \phi_*\tau.\]
    \looseness=-1
    For a set $S$, the unit $\eta_S \colon S \hookrightarrow \gamma S$ maps an element $s\in S$ to the principal filter at $s$,
    \[ \eta_S(s) \coloneqq \set{S_0 \subseteq S | s \in S_0},\]
    while the multiplication $\mu_S \colon \gamma^2 S \to \gamma S$ maps a filter $\Sigma$ on $\gamma S$ to the filter 
    \[ \mu_S(\Sigma) \coloneqq \set{S_0 \subseteq S | \set{\sigma \in \gamma S | S_0 \in \sigma} \in \Sigma}.\]
    Since $\phi_*\tau$ is an ultrafilter if so is $\tau$, the same definitions also determine the \emph{ultrafilter monad} $\beta \colon \Set \to \Set$, mapping each set $S$ to the set $\beta S \subseteq \gamma S$ of ultrafilters on it.
\end{defn}

\begin{rmk}\label{rmk:partial-maps}
    The definition of $\phi_*\tau$ makes sense more generally if $\phi$ is a partial function whose domain is $\tau$-large, and moreover $\phi_*\tau = \phi'_*\tau$ for any other such partial function $\phi'$ which agrees with $\phi$ on a $\tau$-large set. 
\end{rmk}

We recall here the fundamental notion of a \emph{reduced product} of sets, a construction generalizing the more familiar notion of an \emph{ultraproduct} to arbitrary filters. We refer the reader to \cite[\S 6]{changModelTheory1990} for an account of their properties and applications in model theory.

\begin{defn}[Reduced products and ultraproducts]\label{defn:reduced-products-and-ultraproducts}
    Let $\sigma$ be a filter on a set $S$. The \emph{reduced product} of an $S$-family $(X_s)_{s\in S}$ of sets with respect to $\sigma$ is the set
    \[\textstyle \int_{s:\sigma} X_s \coloneqq \colim_{S_0\in\sigma^{\op}}\prod_{s\in S_0}X_s,\]
    whose elements are equivalence classes of tuples $(x_s) \in \prod_{s\in S_0}X_s$ for some $S_0 \in \sigma$, identified if they agree on a $\sigma$-large set. This definition gives rise to the \emph{reduced product functor}:
    \[\begin{tikzcd}[column sep = 50pt]
    	{\Set^S} & {\Set .}
    	\arrow["{\int_{s:\sigma}-}", from=1-1, to=1-2]
    \end{tikzcd}\]
    We denote the elements of $\int_{s:\sigma} X_s$ by $(x_s)_{s:\sigma}$, and we refer to them as \emph{$\sigma$-families}. In the case where all sets $X_s$ coincide with a single set $X$, we also speak of the \emph{reduced power} of $X$ with respect to $\sigma$, which we denote by $X^\sigma$. If $\sigma$ is an ultrafilter, we speak of \emph{ultraproducts} and \emph{ultrapowers} in place of reduced products and reduced powers, respectively.
\end{defn}

\begin{rmk}\label{rmk:reduced-products-and-ultraproducts}
    More liberally, we will also speak of reduced products of, and filter-indexed families in, proper \emph{classes}. Note also that the reduced product $\int_{s:\sigma} X_s$ depends, up to a canonical bijection, only on the $\sigma$-family $(X_s)_{s:\sigma}$; therefore, we will often consider reduced products of filter-indexed families of sets.
\end{rmk}

We conclude this subsection by recalling the \emph{dependent tensor product} of filters, first introduced in \cite[Def.\ 1.8]{katetovProductsFilters1968}, which combines a family of filters with respect to a filter on the index set. This construction extends the familiar tensor product of filters dating back at least to \cite[Thm.\ 1.10]{frayneReducedDirectProducts1962}.

\begin{cons}[Dependent tensor products]\label{cons:dependent-tensor-products}
\looseness=-1
    For a filter $\sigma$ on a set $S$ and an $S$-family $(\tau_s)_{s\in S}$ of filters, each on a set $T_s$, we denote by $\sigma \cdot \tau_s$ the filter on the disjoint union $S \cdot T_s$ given by the subsets $S \cdot T_{s,0}$ satisfying $T_{s,0} \in \tau_s$ for $s :\sigma$. Note that $\sigma \cdot \tau_s$ is an ultrafilter if so are $\sigma$ and each $\tau_s$, and moreover $\sigma \cdot \tau_s = \sigma \cdot \tau'_s$ if $\tau_s = \tau'_s$ for $s:\sigma$. 

    In the case where all filters $\tau_s$ coincide with a single filter $\tau$ on a set $T$, the filter $\sigma\cdot \tau$ is the \emph{tensor product} of $\sigma$ and $\tau$, which we can describe as the filter on the set $S \times T$ consisting of those subsets $R$ satisfying 
    \[ \set{ s\in S | \set{ t \in T | (s,t) \in R } \in \tau } \in \sigma. \]
\end{cons}

\subsection{The category of filters}\label{ssec:category-of-filters}In this subsection, we collect a few properties of the \emph{category of filters}, first defined in \cite{koubekCategoryFilters1970} and studied extensively in \cite{blassTwoClosedCategories1977}; its full subcategory on ultrafilters is also central in Blass's PhD thesis \cite{blassOrderingsUltrafilters1970}. 

\begin{defn}[$\FF$ and $\UF$]
    The \emph{category of filters} $\FF$ is defined as follows: 
    \begin{enumerate}
        \item objects are pairs $(S, \sigma)$ of a set $S$ together with a filter $\sigma$ on $S$;
        \item arrows $(T, \tau) \to (S, \sigma)$ are equivalence classes of partial functions $\phi \colon T \rightharpoonup S$ with $\dom \phi \in\tau$ and such that $\sigma \subseteq \phi_*\tau$, where two such maps are identified if they agree on a $\tau$-large subset.
    \end{enumerate}
    The \emph{category of ultrafilters} $\UF$ is the full subcategory of $\FF$ spanned by ultrafilters.
\end{defn}

Note that, for a filter $(S,\sigma)$ and some $S_0\in\sigma$, the inclusion $S_0 \hookrightarrow S$ realizes an isomorphism $(S_0,\sigma_{\upharpoonright S_0})\cong (S,\sigma)$ in $\FF$: in other words, the carrier set of a filter is not invariant under isomorphism in $\FF$, so that we will typically omit it. Moreover, we will not distinguish notationally between a function and the arrow it defines up to equivalence; note that each arrow $\tau \to \sigma$ is represented up to equivalence by a total function $T \to S$ unless $S$ is empty, in which case $\tau$ is improper and there is a unique such arrow.

\begin{nota}\label{nota:def-deltaS}
    For a set $S$, we denote by $\delta_S$ the filter $\set{S}$. The association $S \mapsto \delta_S$ realizes an embedding $\Set \hookrightarrow \FF$, which we leave implicit when unambiguous.
\end{nota}

Unlike $\UF$, the category of filters $\FF$ has a rich categorical structure: for instance, as shown in \cite[Thm.\ 9]{blassTwoClosedCategories1977}, it admits all finite limits, whose carrier sets are computed as in $\Set$. We also recall the following facts from \cite[Thm.\ 6]{blassTwoClosedCategories1977}.
\begin{enumerate}
    \item\looseness=-1 An $\FF$-arrow $\phi \colon \tau \to \sigma$ is a monomorphism if and only if it is injective on a $\tau$-large subset; in particular, we will always identify a filter $\tau$ refining $\sigma$ with the canonical $\FF$-arrow $\tau \hookrightarrow \sigma$ represented by the identity map.
    \item  An $\FF$-arrow $\phi \colon \tau \to \sigma$ is a cover if and only if $\phi(T_0) \in \sigma$ for each $T_0 \in \tau$, i.e.\ if and only if $\sigma = \phi_*\tau$. 
    \item $\FF$ is \emph{balanced}, i.e.\ an arrow which is both monic and epic is an isomorphism.
\end{enumerate}

\looseness=-1
Note also that epimorphisms in $\FF$ --- henceforth called \emph{covers} --- are effective and pullback-stable (cf.\ \cite[Lem.\ 1.2]{moerdijkModelIntuitionisticNonstandard1995}); therefore, $\FF$ is a regular category where all epimorphisms are regular. Every $\FF$-arrow whose domain is proper and whose codomain is an ultrafilter is a cover, so in particular all arrows in $\UF$ are covers.

\begin{prop}[Coproducts and extensivity]\label{prop:F-is-lextensive}
\looseness=-1
    The category $\FF$ admits small coproducts; moreover, it is (finitary) \emph{extensive}, meaning that finite coproducts are disjoint and stable under pullback.
\begin{proof}
    The coproduct of a family $(\sigma_i)_{i\in I}$ of filters is given by the dependent tensor product $I \cdot \sigma_i$ (cf.\ \Cref{cons:dependent-tensor-products}); the fact that coproducts are disjoint then follows using the explicit description of pullbacks. Moreover, it is clear that the (improper) filter $0$ on the empty set is a strict initial object. Finally, for $\phi\colon \sigma\to 2$, seeing $2$ as $\set{0,1}$, one can easily compute that $\sigma = \sigma_0+\sigma_1$, where $\sigma_i$ is the pullback of the canonical inclusion $\set{i}\mono2$ along $\phi$, which gives an equivalence $\FF \times \FF \simeq \FF/2$. We can deduce the general case $\FF/\tau_0 \times \FF/\tau_1 \simeq \FF/(\tau_0+\tau_1)$ for any pair of filters $\tau_0,\tau_1$, proving extensivity, by composing any $\sigma\to \tau_0+\tau_1$ with the canonical map $\tau_0+\tau_1\to2$. 
\end{proof}
\end{prop}

\begin{rmk}\label{prop:F-infextensive}
    Note that $\FF$ is not \emph{infinitary} extensive. Indeed, for $\sigma$ the Fréchet filter on $\mathbb N$, the pullback of any $\set{n}\mono\mathbb N$ along $\sigma\mono\mathbb N$ is the initial filter $0$. Since $\mathbb N \cong \mathbb N \cdot \set{n}$, the arrows $\set{n}\mono\mathbb N$ can be seen as coprojections into a countable coproduct: infinitary extensivity would then force $\sigma \cong 0$.
\end{rmk}

It is clear that the above properties carry over to the slice categories $\FF/\sigma$ for any filter $\sigma$. In particular, as described in \cite[Ex.\ C2.1.12(b)]{johnstoneSketchesElephantTopos2002}, we can consider the \emph{regular topology} on $\FF/\sigma$, which here amounts to the Grothendieck topology generated by covers. This yields the following explicit description of the sheaf condition (cf.\ \cite[\S I.4]{barrExactCategories1971}), which will play a crucial role in \Cref{sec:tautness}.

\begin{prop}[Sheaves on $\FF/\sigma$]\label{prop:sheaves-on-F-over-sigma}
    A presheaf $\mathrm H \colon (\FF/\sigma)^{\op} \to \Set$ is a sheaf for the regular topology on $\FF/\sigma$ if and only if for every cover $\theta \colon \psi \epi \phi$ over $\sigma$, with kernel pair   
    \[\begin{tikzcd}
    	{\kappa\times_{\lambda}\kappa} & {\kappa} & \lambda \\
    	& \sigma
    	\arrow["{\pi_1}", shift left, from=1-1, to=1-2]
    	\arrow["{\pi_2}"', shift right, from=1-1, to=1-2]
    	\arrow["{\psi\times_\phi\psi}"', curve={height=12pt}, from=1-1, to=2-2]
    	\arrow["\theta", two heads, from=1-2, to=1-3]
    	\arrow["\psi", from=1-2, to=2-2]
    	\arrow["\phi", curve={height=-12pt}, from=1-3, to=2-2]
    \end{tikzcd},\]
    the diagram
    \[\begin{tikzcd}
    	{\mathrm H\phi} & {\mathrm H\psi} & {\mathrm H(\psi\times_\phi \psi)}
    	\arrow["{\mathrm H\theta}", from=1-1, to=1-2]
    	\arrow["{\mathrm H\pi_1}", shift left, from=1-2, to=1-3]
    	\arrow["{\mathrm H\pi_2}"', shift right, from=1-2, to=1-3]
    \end{tikzcd}\]
    is an equalizer in $\Set$.
\end{prop}

\begin{rmk}[The coherent topology]\label{rmk:F-is-coherent}
    The above discussion also entails that $\FF$ is a coherent category (see also \cite[\S A.5]{pittsConceptualCompletenessFirstorder1989}). This means that we can also endow $\FF/\sigma$ with the \emph{coherent topology}, the combination of the regular topology with the extensive topology described, e.g., in \cite[Ex.\ C2.1.12(d)]{johnstoneSketchesElephantTopos2002}; this topology on $\FF$ has already been considered in \cite{moerdijkModelIntuitionisticNonstandard1995,palmgrenRealNumbersTopos2001,awodeyUltrasheavesDoubleNegation2004}.
    We note here that a presheaf $\mathrm H \colon (\FF/\sigma)^{\op} \to \Set$ is a sheaf for the coherent topology if and only if it satisfies the condition in \Cref{prop:sheaves-on-F-over-sigma} and moreover it maps finite coproducts to products.
\end{rmk}

Using the category of filters $\FF$, we can describe reduced products in a more intrinsic way; let us start with reduced powers.

\begin{cons}[Reduced powers]\label{cons:reduced-powers}
\looseness=-1
For a filter $\sigma$ and a set $X$, the reduced power $X^\sigma$ is in bijection with the hom-set $\FF( \sigma, X)$. For an $\FF$-arrow $\phi \colon \tau \to \sigma$, the precomposition map $\phi^* \colon \FF(\sigma, X) \to \FF(\tau,X)$ then corresponds to the function $(x_{s})_{s:\sigma} \mapsto (x_{\phi t})_{t:\tau}$, so that the reduced power functor $X^{-}\colon \FF^{\op} \to \Set$ is isomorphic to the presheaf represented by $\delta_X$.
\end{cons}

In order to generalize the above to reduced products, we have to work inside \emph{slices} of $\FF$.

\begin{cons}[Reduced products]\label{cons:reduced-products}
For an $S$-family $(X_s)_{s\in S}$ of sets and a filter $\sigma$ on $S$, the reduced product $\int_{s:\sigma} X_s$ is in bijection with the hom-set $\FF/S (\sigma, S\cdot X_s)$, i.e.\ the set of sections in $\FF$ of the evident projection $S\cdot X_s \to S$ `with respect to $\sigma$':
\[\begin{tikzcd}[sep =20pt]
	& {S\cdot X_s} \\
	\sigma & S
	\arrow[from=1-2, to=2-2]
	\arrow[curve={height=-12pt}, dashed, from=2-1, to=1-2]
	\arrow[hook, from=2-1, to=2-2]
\end{tikzcd}\] 
As above, an $\FF$-arrow $\phi \colon \tau \to \sigma$, represented by a total map $\phi \colon T\to S$, induces a precomposition map $\phi^* \colon \FF/S(\sigma, S\cdot X_s) \to \FF/S(\tau, S\cdot X_s)$; since the square 
\[\begin{tikzcd}[sep =20pt]
	{T \cdot X_{\phi t}} & {S\cdot X_s} \\
	T & S
	\arrow[from=1-1, to=1-2]
	\arrow[from=1-1, to=2-1]
	\arrow["\lrcorner"{anchor=center, pos=0.125}, draw=none, from=1-1, to=2-2]
	\arrow[from=1-2, to=2-2]
	\arrow["\phi"', from=2-1, to=2-2]
\end{tikzcd}\]
is a pullback, the set $\FF/S(\tau, S\cdot X_s)$ is in bijection with $\FF/T(\tau, T\cdot X_{\phi t})$ and hence with the reduced product $\int_{t:\tau} X_{\phi t}$, so that $\phi^*$ corresponds to the function $(x_{s})_{s:\sigma} \mapsto (x_{\phi t})_{t:\tau}$.
\end{cons}

The dependent tensor product of \Cref{cons:dependent-tensor-products} can also be analyzed in terms of $\FF$-arrows as follows.

\begin{cons}[Functoriality of dependent tensor products]\label{rmk:functoriality-dependent-tensor-products}
The dependent tensor product is functorial in the following sense. Let $\sigma$ be a filter and $(\tau_s)_{s:\sigma}$ be a $\sigma$-family of filters.
\begin{enumerate}
    \item For a $\sigma$-family of $\FF$-arrows $(\psi_s \colon \kappa_s \to \tau_s)_{s :\sigma}$, the map $(s, k) \mapsto (s, \psi_s k)$ defines an $\FF$-arrow $\sigma \cdot \kappa_s \to \sigma \cdot \tau_s$, which we denote by $\operatorname{id}_\sigma \cdot\,  \psi_s$.
    \item For an $\FF$-arrow $\phi \colon \lambda \to \sigma$, the map $(\ell,t) \mapsto (\phi \ell, t)$ defines an $\FF$-arrow $\lambda \cdot \tau_{\phi \ell} \to \sigma \cdot \tau_s$, which we denote by $\phi \cdot \operatorname{id}_{\tau_s}$.
\end{enumerate}
These assignments satisfy evident coherence axioms that we omit here. 
\end{cons}

\begin{nota}[Dependent tensor product over a filter]\label{rmk:dependent-tensor-products-on-slices}
    For a $\lambda$-family of $\FF$-arrows $(\phi_\ell \colon \tau_\ell \to \sigma)_{\ell :\lambda}$ over a common filter $\sigma$, we denote by $\lambda \cdot \phi_\ell$ the $\FF$-arrow $\lambda \cdot \tau_\ell \to \sigma$ obtained as the composite of $\operatorname{id}_\lambda \cdot\,  \phi_\ell \colon \lambda \cdot \tau_{\ell} \to \lambda \cdot \sigma$ with the projection $\lambda\cdot\sigma\to\sigma$, i.e.\ the map defined by $(\ell, t) \mapsto \phi_\ell t$.
\end{nota}

The classical fact that every filter coincides with the intersection of all the ultrafilters refining it here takes the form of the following lemma, which can be seen as a categorical formulation of the ultrafilter principle.

\begin{lem}[Ultrafilter principle]\label{lem:cover-by-uf}
    Let $L$ be a set indexing the ultrafilters $(\iota_\ell\colon \sigma_{\ell} \hookrightarrow \sigma)_{\ell \in L}$ refining $\sigma$. Then, the $\FF$-arrow $L\cdot\iota_\ell\colon L\cdot\sigma_\ell\to \sigma$ is a cover.
\begin{proof}
    Straightforward (cf.\ also \cite[Lem.\ 3.4]{awodeyUltrasheavesDoubleNegation2004}).
\end{proof}
\end{lem}

\begin{cons}[Monoidal structure on $\FF$]\label{rmk:monoidal-structure-on-F}
    The tensor product of filters defines a monoidal structure on $\FF$ which, as explained in \cite[\S 4]{blassTwoClosedCategories1977}, is \emph{right closed}: 
    given two filters $(T,\tau)$ and $(L,\lambda)$, the internal hom-object $[\tau, \lambda]$ is the filter on the set $\FF(\tau,L)$ generated by the subsets $\FF(\tau,L_0)$ for $L_0 \in \lambda$. This means that for any filter $\sigma$,
    \[\FF(\sigma,[\tau,\lambda])\cong\FF(\sigma\cdot\tau,\lambda),\]
    and in particular $[\tau, \delta_X] = \delta_{\FF(\tau, X)}$ for any set $X$, so that
    \[ \FF(\sigma, \FF(\tau, X)) \cong \FF(\sigma \cdot \tau , X). \]
    In other words, we get the following `power rule' for reduced powers: 
    \[ (X^\tau)^\sigma \cong X^{\sigma \cdot \tau};\]
    an analogous bijection, reminiscent of the \emph{Fubini rule} for double integrals, can be deduced more generally for reduced products:
    \[\textstyle \int_{s:\sigma} \int_{t:\tau_s} X_{s,t} \cong \int_{(s,t) : \sigma\cdot \tau_s} X_{s,t}.\]
\end{cons}

We conclude this subsection with the following observation, anticipating the constructions of the next one. In a sense that we will make precise later, in \Cref{rmk:Fop-is-categorical}, dependent tensor products endow the category $\FF^{\op}$ with reduced products. As in the case of $\Set$ (\Cref{defn:reduced-products-and-ultraproducts}), we can show that these reduced products are computed as directed colimits of products.

\begin{prop}[Dependent tensor products as reduced products]\label{prop:tensor-is-categorical}
    Let $\sigma$ be a filter on a set $S$ and $(\tau_s)_{s\in S}$ be an $S$-family of filters. Then:
    \[\sigma\cdot\tau_s\cong\lim_{S_0\in\sigma}(S_0\cdot\tau_s).\]
\end{prop}
\begin{proof}
First note that the family of $\FF$-arrows $(\sigma\cdot\tau_s\mono S_0\cdot\tau_s)_{S_0 \in \sigma}$ is a cone over the diagram $\sigma \to \FF$ mapping $S_0\mapsto S_0\cdot\tau_s$. To see that it is a limit, suppose that $\lambda$ is another cone in $\FF$, witnessed by $\FF$-arrows $(\phi_{S_0} \colon \lambda\to S_0\cdot\tau_s)_{S_0\in\sigma}$. Since each $\sigma\cdot\tau_s\mono S_0\cdot\tau_s$ is a monomorphism, it suffices to show that each arrow $\phi_{S_0}$ factors through $\sigma\cdot\tau_s$. This follows from the fact that the diagram 
    \[\begin{tikzcd}
    	{\sigma\cdot\tau_s} & {S_0\cdot\tau_s} \\
    	\sigma & S_0
    	\arrow[hook, from=1-1, to=1-2]
    	\arrow[from=1-1, to=2-1]
    	\arrow["\lrcorner"{anchor=center, pos=0.125}, draw=none, from=1-1, to=2-2]
    	\arrow[from=1-2, to=2-2]
    	\arrow[hook, from=2-1, to=2-2]
    \end{tikzcd}\]
can be seen to be a pullback in $\FF$, together with the fact that the image of $\lambda\xrightarrow{\phi_{S_0}} S_0\cdot\tau_s \to S_0$ contains all $\sigma$-large sets and hence it factors through $\sigma \hookrightarrow S_0$.
\end{proof}

\subsection{Categories with reduced products}\label{ssec:completions}
\looseness=-1
As we have seen in \Cref{def:filter-monads}, mapping each set to its set of filters determines a monad $\gamma$ on $\Set$, and similarly ultrafilters determine the monad $\beta$. In this subsection, we will see that these constructions can be suitably \emph{categorified}, in that there exist two \emph{pseudomonads} $\bbgamma$ and $\bbbeta$ on $\CAT$ which we can think of as analogues on categories of the filter and ultrafilter monads. Intuitively, these will be defined by `freely endowing' a category $\mathcal C$ with reduced products, or ultraproducts, represented by (ultra)filter-indexed families of objects in $\mathcal C$: for this reason, we name them the \emph{reduced-product completion} and the \emph{ultraproduct completion}, respectively. For the relevant notions about pseudomonads we refer the reader to, e.g., \cite{lackCoherentApproachPseudomonads2000}.

In the case of ultrafilters, the pseudomonad $\bbbeta$ was first introduced in \cite{rosoliniUltracompletions2024}, building on Garner's work \cite{garnerUltrafiltersFiniteCoproducts2020}, to study Makkai's \emph{ultracategories} \cite{makkaiStoneDualityFirst1987}. Indeed, as solidified by the results of \cite{hamadUltracategoriesColaxAlgebras2025}, its \emph{pseudoalgebras} are a reasonable axiomatization of the intuitive notion of an ultracategory as a category equipped with ultraproducts; we refer the reader to \cite[\S 2]{tarantinoUltracategoriesKanExtensions2025} for a discussion of the subject. In the same way, pseudoalgebras for $\bbgamma$ --- henceforth simply called \emph{algebras} --- can be thought of as categories equipped with reduced products.

\begin{cons}[Completions]\label{cons:reduced-product-completion}
For a category $\mathcal C$, we denote by $\bbgamma \mathcal C$ the category defined as follows:
\begin{enumerate}
    \item objects are pairs of a filter $\sigma$ and a $\sigma$-family $(c_s)_{s:\sigma}$ of objects in $\mathcal C$;
    \item arrows $(c_{s})_{s:\sigma} \to (d_t)_{t:\tau}$ are pairs of an $\FF$-arrow $\phi \colon \tau \to \sigma$ and a $\tau$-family of $\mathcal C$-arrows $(f_t \colon c_{\phi t} \to d_t)_{t: \tau}$.
\end{enumerate}
We also denote by $\bbbeta \mathcal C$ its full subcategory spanned by ultrafilter-indexed families.

\looseness=-1
Both constructions give rise to pseudomonads on $\CAT$ which we refer to as the \emph{reduced-product completion} $\bbgamma$ and the \emph{ultraproduct completion} $\bbbeta$. For a category $\mathcal C$, 
\begin{enumerate}
    \item the unit $\bbeta_\mathcal C \colon \mathcal C \to \bbgamma \mathcal C$ maps an object $c$ in $\mathcal C$ to the $1$-family $(c)_{*:1}$,
    \item and the multiplication $\bbmu_\mathcal C \colon \bbgamma^2 \mathcal C \to \bbgamma \mathcal C$ maps a $\sigma$-family of filter-families $( (c_{s,t})_{t:\tau_s})_{s:\sigma}$ in $\mathcal C$ to the $(\sigma\cdot\tau_s)$-family $(c_{s,t})_{(s,t):\sigma\cdot\tau_s}$;
\end{enumerate}
both assignments are evidently functorial and restrict to a unit and multiplication for $\bbbeta$. In particular, the inclusions $\bbbeta \mathcal C \hookrightarrow \bbgamma \mathcal C$ determine a pseudomonad morphism $\bbbeta \to \bbgamma$ in the sense of, e.g., \cite[Def.\ 2.1]{gambinoFormalTheoryPseudomonads2021}.

Note also that $\bbgamma$ recovers the category of filters as $\FF\simeq (\bbgamma 1)^{\op}$, and similarly $\UF \simeq (\bbbeta 1)^{\op}$.
\end{cons}

Let us describe two examples of $\bbgamma$-algebras which will play a central role in the sequel; throughout, we will speak of the \emph{reduced product} of a $\sigma$-family of objects $(c_s)_{s:\sigma}$ in a $\bbgamma$-algebra $\mathcal C$ as the image of $(c_s)_{s:\sigma}$ under the algebra map $\bbgamma \mathcal C \to \mathcal C$.

\begin{ex}[$\FF^{\op}$ and $(\FF/\sigma)^{\op}$]\label{ex:slices-of-F-as-bbgamma-algebras}
    Through the equivalence with $\bbgamma 1$, the category $\FF^{\op}$ carries the structure of a $\bbgamma$-algebra, where the reduced product of a $\sigma$-family of filters $(\tau_s)_{s:\sigma}$ is given by their dependent tensor product $\sigma \cdot \tau_s$; functoriality of the algebra map $\bbgamma(\FF^{\op}) \to \FF^{\op}$ recovers the functoriality of dependent tensor products discussed in \Cref{rmk:functoriality-dependent-tensor-products}.
    More generally, for any filter $\sigma$, the category $(\FF/\sigma)^{\op}$ also carries the structure of a $\bbgamma$-algebra, where the reduced product of a $\lambda$-family of $\FF$-arrows $(\phi_\ell \colon \tau_\ell \to \sigma)_{\ell:\lambda}$ is given by the $\FF$-arrow $\lambda \cdot \phi_\ell \colon \lambda \cdot \tau_\ell \to \sigma$ (\Cref{rmk:dependent-tensor-products-on-slices}).
\end{ex}

\begin{ex}[$\Set$]\label{ex:set-as-a-bbgamma-algebra}
The category $\Set$ carries a natural structure of a $\bbgamma$-algebra, where the reduced product of a $\sigma$-family of sets $(X_s)_{s:\sigma}$ is given by their actual reduced product $\int_{s:\sigma} X_s$ in the sense of \Cref{defn:reduced-products-and-ultraproducts} (see also \Cref{rmk:reduced-products-and-ultraproducts}).
\end{ex}

\begin{ex}[Categorical structures]\label{ex:horn-theories}\label{ex:categorical-bbgamma-alg}
    Abstracting from the case of $\Set$, we say that a $\bbgamma$-algebra structure on a category $\mathcal C$ is \emph{categorical} if the reduced product of a $\sigma$-family of objects $(c_s)_{s:\sigma}$ is given by the usual formula of a directed colimit of products in $\mathcal C$:
    \[\textstyle \colim_{S_0\in\sigma^{\op}}\prod_{s\in S_0} c_s.\]
    This definition endows a category $\mathcal C$ admitting small products and directed colimits with a $\bbgamma$-algebra structure, provided that directed colimits \emph{distribute} over small products in $\mathcal C$.\footnote{The distributivity property is needed for the \emph{multiplicator} of the algebra, corresponding to the Fubini-like bijection of \Cref{rmk:monoidal-structure-on-F}.} In other words, every \emph{predefinable category} in the sense of \cite[Def.\ 2.2.1]{kuberDefinableCategories2018} admits a categorical $\bbgamma$-algebra structure: in particular, so does the category of models of a \emph{regular theory}, where reduced products are computed as in $\Set$.
\end{ex}

\begin{rmk}\label{rmk:Fop-is-categorical}
    By \Cref{prop:tensor-is-categorical}, the $\bbgamma$-algebra structure on $\FF^{\op}$ described in \Cref{ex:slices-of-F-as-bbgamma-algebras} is categorical. Since $\FF^{\op}$ does not admit all directed colimits (cf.\ \cite[Ex.\ 12]{blassTwoClosedCategories1977}), it is a categorical $\bbgamma$-algebra which is not predefinable.
\end{rmk}

Moving on to $\bbbeta$-algebras, we first observe that any $\bbgamma$-algebra canonically yields a $\bbbeta$-algebra by virtue of the morphism of pseudomonads $\bbbeta \hookrightarrow \bbgamma$; see, e.g., \cite[Lem.\ 3.1]{gambinoFormalTheoryPseudomonads2021}. It is also evident that $\UF^{\op}$ (and more generally $(\UF/\sigma)^{\op}$ for any ultrafilter $\sigma$) carries a natural $\bbbeta$-algebra structure, paralleling \Cref{ex:slices-of-F-as-bbgamma-algebras}, which does not arise in this way from a $\bbgamma$-algebra structure. Below, we describe another such example of interest, motivating the aforementioned theory of ultracategories.

\begin{ex}[Models of coherent theories]\label{ex:first-order-theories}
    Let $\mathbb T$ be a coherent theory. The category $\operatorname{Mod}(\mathbb T)$ of models of $\mathbb T$ and their homomorphisms carries a natural structure of a $\bbbeta$-algebra, where ultraproducts are computed as in $\Set$: this follows since, by \L o\'s's theorem, the class of models of a coherent theory is closed under ultraproducts; however, it need not be closed under \emph{reduced} products, so that this $\bbbeta$-algebra structure does not generally arise from a more general $\bbgamma$-algebra structure. 
\end{ex}

We conclude this subsection with a family of examples of $\bbgamma$-algebras and $\bbbeta$-algebras of a quite different nature from those above.

\begin{ex}[Categories of metric spaces]\label{ex:metric-spaces}
    For $k \in (0,\infty]$, let $\mathsf{Met}_k$ be the category of \emph{$k$-bounded metric spaces}, i.e.\ metric spaces\footnote{Crucially, allowing infinite distances if $k = \infty$.} $(X,d_X)$ such that $d_X(x,y) \leq k$ for all $x,y\in X$, and \emph{1-Lipschitz maps}, i.e.\ functions $f \colon X \to Y$ such that $d_Y(f(x), f(y)) \leq d_X(x,y)$ for all $x,y\in X$; let $\mathsf{CompMet}_k$ be the full subcategory of $\mathsf{Met}_k$ spanned by the \emph{complete} metric spaces.

    For a $\sigma$-family $(X_s)_{s:\sigma}$ of $k$-bounded metric spaces, we define its \emph{reduced product} $\int_{s:\sigma} X_s$ as the $k$-bounded metric space defined as follows. Points of $\int_{s:\sigma} X_s$ are equivalence classes of tuples $(x_s)\in \prod_{s\in S_0}X_s$ for some $S_0 \in \sigma$, where $(x_s)$ and $(y_s)$ are identified if, for every $\epsilon>0$, it holds that $\forall_{s:\sigma}d_{X_s}(x_s,y_s)<\epsilon$. Then, the distance on $\int_{s:\sigma} X_s$ is defined by
    \[ \textstyle d_\sigma((x_s),(y_s)) \coloneqq \limsup_{s:\sigma} d_{X_s} (x_s,y_s) = \inf_{S_0 \in \sigma} \sup_{s\in S_0} d_{X_s} (x_s,y_s).\]
    It is straightforward to see that this indeed defines a metric space, which moreover is complete if so is each $X_s$. This operation endows both $\mathsf{Met}_k$ and $\mathsf{CompMet}_k$ with a $\bbgamma$-algebra structure --- which can be shown to be categorical, cf.\ \cite[Prop.\ 2.1]{hamadBundlesMetricStructures2026} --- and hence with a $\bbbeta$-algebra structure; see also \cite[Thm.\ 8.2]{hamadUltracategoriesColaxAlgebras2025}, \cite{henryUltracategoriesUltracategoriesMetric2026}, and \cite{yaacovModelTheoryMetric2008,ghasemiReducedProductsMetric2016}. 
\end{ex}

\subsection{Vf-categories and vu-categories}\label{ssec:vf-cats-vu-cats}In this subsection, we introduce \emph{vf-categories}, the main object of study in this paper, which generalize the notion of a category endowed with reduced products. The definition is inspired by the notion of a \emph{vu-category} introduced in \cite{saadiaExtendingConceptualCompleteness2026,hamadGeneralisedUltracategoriesConceptual2025,vangoolToposesEnoughPoints2026}, which we also recall below. As already discussed in the introduction, we can think of vf-categories and vu-categories as a categorification of the convergence relations induced by topological spaces, between filters and points and between ultrafilters and points, respectively.

\begin{defn}[Vf-categories]\label{defn:vf-cats}
    A \emph{vf-category} $\X$ is defined by the following data:
    \begin{enumerate}
        \item a class $\X_0$, whose elements are called \emph{points};
        \item a functor $\X \colon \X_0 \times \bbgamma \X_0 \to \Set$, where
        \begin{itemize}
        \item[--] the set $\mathscr{X}(p, (q_s)_{s:\sigma})$ is denoted by $\X_\sigma(p,q_s)$;
        \item[--] elements $u \in \X_\sigma(p,q_s)$ are called \emph{vf-arrows} of \emph{type} $\sigma$ and are denoted graphically as $u \colon p \vuto{\sigma} q_s$;
        \item[--] for a vf-arrow $u \colon p \vuto{\sigma} q_s$ and an $\FF$-arrow $\phi \colon \tau \to \sigma$, we denote by $\phi^* u$ the vf-arrow $p \vuto{\tau} q_{\phi t}$ obtained by functoriality with respect to the $\bbgamma \X_0$-arrow $(\phi, \operatorname{id}) \colon (q_s)_{s:\sigma}\to(q_{\phi t})_{t:\tau}$, which we refer to as the \emph{reindexing} of $u$ along $\phi$;
        \end{itemize} 
        \item for each point $p$, an \emph{identity} vf-arrow $\id_p \colon p \vuto{\ptuf} p$;
        \item for each vf-arrow $u \colon p\vuto{\sigma}q_s$ and each $\sigma$-family $(v_s\colon q_s\vuto{\tau_s}r_{s,t})_{s: \sigma}$ of vf-arrows, a \emph{composite} $v_s\circ_\sigma u \colon p \vuto{\sigma\cdot\tau_s} r_{s,t}$, which we also denote graphically as
        \[v_s\circ_\sigma u \colon p \vuto{\sigma} (q_s\vuto{\tau_s} r_{s,t}).\]                
    \end{enumerate}
This data must satisfy the following axioms:
    \begin{enumerate}
        \item unitality and associativity of composition,
            \[u\circ_{\ptuf} \id_p = u, \qquad \id_{q_s}\circ_{\sigma} u = u,\]
            \[ (w_{s,t}\circ_{\tau_s} v_s)\circ_{\sigma} u = w_{s,t}\circ_{\sigma\cdot \tau_s} (v_s\circ_{\sigma} u);\]
        \item compatibility of reindexing with composition,
        \[(\psi_s^*v_s)\circ_{\sigma} u = (\operatorname{id}_{\sigma}\cdot\ \psi_s)^*(v_s\circ_{\sigma} u),\qquad v_{\phi \ell}\circ_{\lambda} (\phi^*u) = (\phi\cdot \operatorname{id}_{\tau_s})^*(v_s\circ_{\sigma} u),\]
         for each $\sigma$-family of $\FF$-arrows $(\psi_s \colon \kappa_s\to\tau_s)_{s:\sigma}$ and each $\FF$-arrow $\phi \colon \lambda \to \sigma$.
    \end{enumerate}

    We will typically write $\X$ in place of $\X_0$.
\end{defn} 

\begin{rmk}
    \looseness=-1
    In the equalities of (i) above, we commit an implicit abuse of notation: for instance, $u$ and $u\circ_{\ptuf} \id_p$ cannot be `equal' as they do not live in the same set. However, they can be identified by reindexing along the canonical isomorphism $1\cdot\sigma\cong\sigma$. Identifications of this kind, already implicit in \cite{saadiaExtendingConceptualCompleteness2026,hamadGeneralisedUltracategoriesConceptual2025,vangoolToposesEnoughPoints2026}, will be made freely in the rest of the paper.
\end{rmk}

\begin{nota}[Specialization and thickening]
    For an $\FF$-arrow $\phi \colon\tau \to \sigma$, we also refer to the reindexing of $u \colon p\vuto{\sigma}q_s$ along $\phi$ as a \emph{specialization} of $u$ if $\phi$ is a monomorphism, and as a \emph{thickening} of $u$ if $\phi$ is a cover. For each point $p$, we denote by $\diag^p_\sigma \colon p \vuto{\sigma} p$ the vf-arrow obtained by reindexing $\id_p \colon p\vuto{\ptuf} p$ along the unique $\FF$-arrow $\sigma \to \ptuf$, which we refer to as the \emph{$\sigma$-diagonal} at $p$.
\end{nota}

\begin{defn}[Vf-functors and vf-transformations]
    \looseness=-1
    For vf-categories $\X$ and $\Y$, a \emph{vf-functor} $F\colon \X\to\Y$ consists of a map on points, associating to each $p\in \X_0$ a point $F(p)\in \Y_0$, and a map on vf-arrows, associating to each $u\colon p\vuto{\sigma}q_s$ in $\X$ a vf-arrow $F(u)\colon F(p)\vuto{\sigma}F(q_s)$ in $\Y$, satisfying the following axioms:
    \[F(\id_p) = \id_{F(p)}, \quad F(v_s)\circ_\sigma F(u) = F(v_s\circ_\sigma u), \quad F(\phi^*u) = \phi^*F(u).\]

    For vf-functors $F,G \colon \X \to \Y$, a \emph{vf-transformation} $\alpha \colon F\Rightarrow G$ consists of a vf-arrow $\alpha_p \colon F(p)\vuto{\ptuf}G(p)$ of type $\ptuf$ in $\Y$ for each $p\in \X_0$, such that 
    \[\alpha_{q_s} \circ_\sigma F(u) = G(u)\circ_1\alpha_p\]
    holds for each vf-arrow $u\colon p\vuto{\sigma}q_s$ in $\X$.
    
    We denote by $\vfCAT$ the 2-category of vf-categories, vf-functors, and vf-transformations. 
\end{defn}

\begin{defn}[Vf-posets]
    A \emph{vf-poset} is a vf-category which is \emph{posetal}, meaning that for each point $p$ and each $\sigma$-family of points $(q_s)_{s:\sigma}$, there is at most one vf-arrow $p \vuto{\sigma}q_s$; if there is one, we also write $p \vuleq{\sigma} q_s$, and we say that $(q_s)_{s:\sigma}$ \emph{converges} to $p$. We denote by $\vfPos$ the full sub-2-category of $\vfCAT$ spanned by vf-posets. Note that we will always assume the class of points of a vf-poset to be, in fact, a set.
\end{defn}

\begin{ex}[Topological spaces]\label{ex:top-sp-as-vf-posets}
    A topological space can be canonically seen as a vf-poset by setting $p \vuleq{\sigma} q_s$ if $q_*\sigma$ converges to $p$, that is, if every neighborhood of $p$ is $q_*\sigma$-large. This definition gives rise to a 2-fully-faithful embedding $\TopSp \mono \vfPos$.
\end{ex}

\begin{ex}[The vf-category of sets]\label{ex:sets}
    The class $\Set$ of sets can be endowed with the structure of a vf-category by setting vf-arrows $X \vuto{s:\sigma} Y_s$ to be functions from $X$ into the reduced product $\int_{s:\sigma} Y_s$. The identity vf-arrow $\id_X \colon X \vuto{\ptuf} X$ is defined by $x \mapsto (x)_{*:\ptuf}$. The reindexing of a vf-arrow $u \colon X  \vuto{s:\sigma} Y_s$ along an $\FF$-arrow $\phi \colon \tau \to \sigma$ is the function 
    \[\begin{tikzcd}
        X & {\int_{s:\sigma} Y_s} & {\int_{t:\tau}Y_{\phi t}}
        \arrow["u", from=1-1, to=1-2]
        \arrow["{\phi^*}", from=1-2, to=1-3]
    \end{tikzcd}.\]
    The composite of a vf-arrow $u \colon X  \vuto{s:\sigma} Y_s$ with a $\sigma$-family of vf-arrows $(v_s \colon Y_s \vuto{t:\tau_s} Z_{s,t})_{s:\sigma}$ is the function 
    \[\begin{tikzcd}[sep =35pt]
        X & {\int_{s:\sigma} Y_s} & {\int_{s:\sigma}\int_{t:\tau_s} Z_{s,t}} &[-30pt] {\int_{(s,t):\sigma\cdot \tau_s}Z_{s,t}}
        \arrow["u", from=1-1, to=1-2]
        \arrow["{\int_{s:\sigma} v_s}", from=1-2, to=1-3]
        \arrow["\cong"{description}, draw=none, from=1-3, to=1-4]
    \end{tikzcd}.\]
    
\end{ex}

\begin{ex}[$\bbgamma$-algebras]\label{ex:bbgamma-algebras}
    Generalizing \Cref{ex:sets}, every $\bbgamma$-algebra can be endowed with the structure of a vf-category by setting vf-arrows $p \vuto{\sigma} q_s$ to be actual arrows in the underlying category from $p$ into the reduced product of the $\sigma$-family $(q_s)_{s:\sigma}$. 
\end{ex}

\begin{ex}[The vf-category of points of a topos]\label{ex:points-of-toposes}
    Generalizing both \Cref{ex:top-sp-as-vf-posets,ex:sets}, the class $\pt(\E)$ of points of a topos $\E$ can be endowed with the structure of a vf-category by setting vf-arrows $p \vuto{s:\sigma} q_s$ to be natural transformations from the inverse image functor $p^*\colon \E\to\Set$ into the \emph{reduced product} $\int_{s:\sigma} q_s^*$, defined as the functor:
    \[\begin{tikzcd}[sep = 45pt]
    	\E & {\Set^S} & \Set
    	\arrow["{\braket{q_s^*}_{s\in S}}", from=1-1, to=1-2]
    	\arrow["{\int_{s:\sigma}-}", from=1-2, to=1-3]
    \end{tikzcd}.\]
    Note that, in general, this functor is not itself (the inverse image of) a point of $\E$. Identities, compositions, and reindexings are defined in the same way as in \Cref{ex:sets}.
\end{ex}

\begin{rmk}[Toposes and geometric logic]
    The vf-category of points of a topos also has a logical reading. Recall that every topos classifies a geometric theory $\mathbb T$ and that, conversely, every geometric theory $\mathbb T$ possesses a classifying topos $\mathcal E_{\mathbb T}$, whose points can be identified with models of $\mathbb T$ (see \cite[Ch.\ X]{maclaneSheavesGeometryLogic1994}). Hence, the vf-category structure on $\pt(\E_{\mathbb T})$ described above can be identified with one on $\operatorname{Mod}(\mathbb T)$: through this identification, vf-arrows $M \vuto{\sigma} N_s$ correspond to \emph{structure homomorphisms} from $M$ into the reduced product of the $\sigma$-family $(N_s)_{s:\sigma}$ of models, even if the latter is not itself a model.
\end{rmk}

Replacing filters with ultrafilters in \Cref{defn:vf-cats}, we recover the notion originally introduced in \cite{saadiaExtendingConceptualCompleteness2026,hamadGeneralisedUltracategoriesConceptual2025,vangoolToposesEnoughPoints2026}.

\begin{defn}[Vu-categories]\label{defn:vu-cats}
    A \emph{vu-category} $\X$ consists of a class $\X_0$ of \emph{points} and a functor $\X \colon \X_0 \times \bbbeta \X_0 \to \Set$ specifying \emph{vu-arrows}, together with identities and compositions satisfying axioms analogous to those defining vf-categories.
\end{defn}

In the same way, we define \emph{vu-functors}, \emph{vu-transformations}, and \emph{vu-posets}, with the obvious notation $\vuCAT$ and $\vuPos$. 

\begin{rmk}
    Note that, since all $\UF$-arrows are covers, all reindexings in a vu-category are thickenings.
\end{rmk}

\begin{cons}[Underlying vu-category of a vf-category]
    By restricting to those vf-arrows whose type is an ultrafilter, every vf-category yields a vu-category, so that there is an evident forgetful 2-functor $\mathsf U \colon \vfCAT \to \vuCAT$, and similarly $\mathsf U \colon \vfPos \to \vuPos$. 
\end{cons}

\begin{ex}
    Through the forgetful 2-functor, \Cref{ex:top-sp-as-vf-posets,ex:sets,ex:bbgamma-algebras,ex:points-of-toposes} all carry over to give examples of vu-categories. Moreover, $\bbbeta$-algebras can be seen as vu-categories exactly as in \Cref{ex:bbgamma-algebras}; cf.\ also  \cite[Ex.\ 5.2(vii)]{saadiaExtendingConceptualCompleteness2026} and \cite[Cor.\ 5.26]{aristoteProfunctorialAlgebras2026}.
\end{ex}

\begin{rmk}[Relational algebras]\label{rmk:relational-algebras}
    That topological spaces embed into vu-categories, with a description analogous to that of \Cref{ex:top-sp-as-vf-posets}, is also directly observed in \cite{saadiaExtendingConceptualCompleteness2026, hamadGeneralisedUltracategoriesConceptual2025, vangoolToposesEnoughPoints2026}. We remark here that the embedding $\TopSp \hookrightarrow \vuCAT$ can also be described in algebraic terms, making use of Barr's characterization of topological spaces as the \emph{relational algebras} for the ultrafilter monad $\beta$ \cite{barrRelationalAlgebras1970}. Indeed, a relational algebra $\xi \colon \beta X \pro X $ --- i.e.\ a set $X$ together with a relation $\xi \subseteq X \times \beta X $ satisfying appropriate axioms --- gives rise to a vu-poset on $X$ simply by setting $p \vuleq{\sigma} q_s$ if $p \ \xi\ q_* \sigma$.

    \looseness=-1
    Perhaps counterintuitively, an analogous phenomenon does not occur in the case of filters and vf-categories: relational algebras for the filter monad $\gamma$ need not give rise to vf-categories with the above definition (see also \cite{dayFilterMonadsContinuous1975}). The reason is that relational $\gamma$-algebras `lack specializations': that is, reindexing is only defined along covers in $\FF$. Ultimately, this is due to the fact that the set $\gamma S$ of filters on any set $S$ carries the refinement order, which determines monomorphisms in $\FF$ and which is trivial between ultrafilters. Relational $\gamma$-algebras are \emph{blind} to this order: likewise, the \emph{convergence spaces} of \cite[Def.\ 5.28]{aristoteProfunctorialAlgebras2026} can be thought of as vf-categorical structures `without specializations', and indeed they include relational $\gamma$-algebras. Instead, we believe that vf-categories coincide with the mathematical entity hinted at in \cite[Ques.\ 5.31]{aristoteProfunctorialAlgebras2026}, which does take the natural order of filters into account.
\end{rmk}

\subsection{Merging vf-arrows}\label{ssec:merging}In this last introductory subsection, we define \emph{merging}, an operation on vf-arrows in a vf-category obtained by combining diagonals and compositions, and we analyze it from an algebraic perspective. This analysis will play a key role in the next section, towards an understanding of \emph{tautness} for vf-categories. As we will show in \Cref{lem:merging-is-inverse-to-specialization}, merging also provides an inverse to specialization along the coprojections into a coproduct of filters. Throughout, fix a vf-category $\X$.

\begin{defn}[Merging]\label{cons:merging}
    For a $\lambda$-family $(u_\ell \colon p \vuto{\sigma_\ell} q_{\ell,s})_{\ell: \lambda}$ of vf-arrows with the same domain, their \emph{merge} is the vf-arrow $p \vuto{\lambda \cdot \sigma_\ell} q_{\ell,s}$ obtained as the composite
    \[ u_\ell \circ_\lambda \diag^p_\lambda \colon p \vuto{\lambda} (p \vuto{\sigma_\ell} q_{\ell,s}),\]
    which we denote by $\bigmrg_{\ell:\lambda} u_\ell$. 
\end{defn}

\begin{nota}[Merging pairs]\label{nota:merging-pairs}
    In particular, merging along the filter $2$ (that is, $\delta_2$) gives us a way to combine two vf-arrows of common domain $u_0\colon p\vuto{\sigma}q_s$ and $u_1\colon p\vuto{\tau}r_t$ into one vf-arrow of type $\sigma + \tau$, which we denote by
    \[u_0\mrg u_1\colon p\vuto{\sigma+\tau}((q_s),(r_t)).\]
\end{nota}

To capture the algebraic content of merging, we introduce here the family of \emph{reindexing functors} associated to any vf-category.

\begin{defn}[Reindexing functors]
    Given a point $p$ and a $\sigma$-family of points $(q_s)_{s:\sigma}$ in $\X$, we define the \emph{reindexing functor} at $p$ and $(q_s)_{s:\sigma}$ as the functor:
    \[ (\FF/\sigma)^{\op} \to \Set \qquad (\phi \colon \tau \to \sigma) \mapsto \X_\tau (p, q_{\phi t}).\]   
\end{defn}

Recall from \Cref{ex:slices-of-F-as-bbgamma-algebras} and \Cref{ex:set-as-a-bbgamma-algebra} that the categories $(\FF/\sigma)^{\op}$ for a filter $\sigma$ and $\Set$ naturally carry $\bbgamma$-algebra structures. As we will show in \Cref{prop:merg-struct}, the merging operation of \Cref{cons:merging} endows each reindexing functor $(\FF/\sigma)^{\op} \to \Set$ of $\X$ with the structure of a \emph{lax morphism} of $\bbgamma$-algebras (see, e.g., \cite[Def.\ 2.4.6]{stepanLaxStructures2category2025}), which we consider here under the name of a \emph{merging structure}. 

\begin{defn}[Merging structures]\label{defn:merging-structures} 
    For a filter $\sigma$ and a functor $\mathrm H \colon (\FF/\sigma)^{\op}\to\Set$, a \emph{merging structure} on $\mathrm H$ is a family of maps
    \[\textstyle m_{\lambda} \colon \int_{\ell:\lambda}\mathrm H(\phi_\ell)\to \mathrm H(\lambda\cdot\phi_\ell),\]
    for each $\lambda$-family $(\phi_\ell \colon \tau_\ell\to\sigma)_{\ell:\lambda}$ in $\FF/\sigma$,
    satisfying the following axioms:
    \begin{enumerate}
        \item \emph{unitality}, 
        \[\begin{tikzcd}[sep = 12pt]
        	{\int_{\ptuf}\mathrm H(\phi)} && {\mathrm H(\ptuf\cdot\phi)} \\
        	& {\mathrm H(\phi)} & 
        	\arrow["\cong"', from=1-1, to=2-2]
        	\arrow["\cong", from=1-3, to=2-2]
        	\arrow["{m_\ptuf}", from=1-1, to=1-3]
        \end{tikzcd}\]
        for each $\FF$-arrow $\phi \colon \tau \to \sigma$;
        \item \emph{associativity},
        \[\begin{tikzcd}[column sep = 32pt]
        	{\int_{\kappa\cdot\lambda_k}\mathrm H(\phi_{k,\ell})} && {\mathrm H((\kappa\cdot\lambda_k)\cdot\phi_{k,\ell})} \\
        	{\int_{\kappa}\int_{\lambda_k}\mathrm H(\phi_{k,\ell})} & {\int_{\kappa}\mathrm H(\lambda_k\cdot\phi_{k,\ell})} & {\mathrm H(\kappa\cdot (\lambda_k\cdot \phi_{k,\ell}))}
        	\arrow["m_{\kappa\cdot\lambda_k}", from=1-1, to=1-3]
        	\arrow["\cong"', from=1-1, to=2-1]
        	\arrow["\int_{\kappa}m_{\lambda_k}"', from=2-1, to=2-2]
        	\arrow["m_\kappa"', from=2-2, to=2-3]
        	\arrow["\cong", from=1-3, to=2-3]
        \end{tikzcd}\]
        for each filter $\kappa$ and each $\kappa$-family of filters $(\lambda_k)_{k:\kappa}$;
        \item \emph{naturality}, 
        \[\begin{tikzcd}
        	{\int_{\lambda}\mathrm H(\phi_\ell)} & {\mathrm H(\lambda\cdot\phi_\ell)} \\
        	{\int_{\lambda}\mathrm H(\psi_\ell)} & {\mathrm H(\lambda\cdot\psi_\ell)}
        	\arrow["m_\lambda", from=1-1, to=1-2]
        	\arrow["{\int_{\lambda}\mathrm H(\theta_\ell)}"', from=1-1, to=2-1]
        	\arrow["{\mathrm H(\operatorname{id}_\lambda\cdot\theta_\ell)}", from=1-2, to=2-2]
        	\arrow["m_\lambda"', from=2-1, to=2-2]
        \end{tikzcd} \quad \mbox{ and }\quad
        \begin{tikzcd}
        	{\int_{\lambda}\mathrm H(\phi_\ell)} & {\mathrm H(\lambda\cdot\phi_\ell)} \\
        	{\int_{\kappa}\mathrm H(\phi_{\alpha k})} & {\mathrm H(\kappa\cdot\phi_{\alpha k})}
        	\arrow["m_\lambda", from=1-1, to=1-2]
        	\arrow["{\alpha^*}"', from=1-1, to=2-1]
        	\arrow["{\mathrm H(\alpha\cdot \operatorname{id}_{\phi_\ell})}", from=1-2, to=2-2]
        	\arrow["m_\kappa"', from=2-1, to=2-2]
        \end{tikzcd}\]
        for each $\lambda$-family of $(\FF/\sigma)$-arrows $(\theta_\ell\colon \psi_\ell\to\phi_\ell)_{\ell:\lambda}$ and each $\FF$-arrow $\alpha\colon \kappa\to\lambda$.
    \end{enumerate}
With a slight abuse of language, we will also refer to the functor $\mathrm H \colon (\FF/\sigma)^{\op}\to\Set$ itself as a \emph{merging structure}.
\end{defn}

\begin{prop}[Merging structures on reindexing functors]\label{prop:merg-struct}
    For a point $p$ and a $\sigma$-family of points $(q_s)_{s:\sigma}$ of $\X$, the reindexing functor 
    \[ (\FF/\sigma)^{\op} \to \Set \qquad 
        (\phi \colon \tau\to\sigma) \mapsto \X_{\tau}(p,q_{\phi t}) \]
    carries a merging structure given by the maps
    \[ \textstyle\int_{\ell:\lambda}\X_{\tau_\ell}(p,q_{\phi_\ell t})\to \X_{\lambda\cdot\tau_\ell}(p,q_{\phi_\ell t}) \qquad (u_\ell)_{\ell:\lambda} \mapsto \bigmrg_{\ell:\lambda}u_\ell,\]
    for each $\lambda$-family $(\phi_\ell \colon \tau_\ell \to \sigma)_{\ell:\lambda}$ in $\FF/\sigma$.
\end{prop}
\begin{proof}
Axiom (i) of \Cref{defn:merging-structures} follows from unitality of identity vf-arrows. Axiom (ii) follows from associativity of composition of vf-arrows together with the fact that  $\diag^p_{\kappa\cdot\lambda_k} = \diag^p_{\lambda_k}\circ_{\kappa}\diag^p_\kappa$. Axiom (iii) follows from the two conditions expressing compatibility of reindexing with composition.
\end{proof}

\begin{rmk}[Right ultrafunctors]
    In the context of ultracategories, lax morphisms of $\bbbeta$-algebras coincide with the \emph{right ultrafunctors} of \cite[Def.\ 8.1.1]{lurieUltracategories2018} (cf.\ \cite[Rem.\ 6.3]{tarantinoUltracategoriesKanExtensions2025}). The obvious analogue of \Cref{prop:merg-struct} for vu-categories then expresses that the reindexing functors $(\UF/\sigma)^{\op} \to \Set$ of a vu-category are right ultrafunctors.
\end{rmk}

\begin{prop}\label{lem:merg-struct-preserve-coproducts}
    Let $\mathrm H \colon (\FF/\sigma)^{\op}\to\Set$ be a functor carrying a merging structure. Then, $\mathrm H$ maps small coproducts in $\FF/\sigma$ to small products in $\Set$.
\end{prop}
\begin{proof}
    Fix a family $(\phi_i \colon \tau_i\to\sigma)_{i\in I}$ of $\FF$-arrows. We denote by $\tau$ the coproduct $I \cdot \tau_i$ in $\FF$, so that the coproduct of the family $(\phi_i)_{i\in I}$ in $\FF/\sigma$ can be identified with the canonical $\FF$-arrow $\phi \colon \tau \to \sigma$; the colimiting cocone is determined by the family of coprojections $\iota_i \colon \tau_i \mono \tau$ in $\FF$, seen as $(\FF/\sigma)$-arrows $\phi_i \to \phi$.
    
    We then consider the component 
    \[ m_I \colon \textstyle\prod_{i\in I} \mathrm H\phi_i \to \mathrm H\phi \]
    of the merging structure and we show that this map is a bijection, with inverse given by the map $\mathrm H\phi \to \prod_{i\in I} \mathrm H\phi_i $ canonically induced by the maps $(\mathrm H \iota_i \colon \mathrm H\phi \to \mathrm H\phi_i)_{i\in I}$.

    Fix $i_0\in I$ and identify it with the inclusion $\set{ i_0} \hookrightarrow I$. By unitality and the second naturality axiom of the merging structure, the diagram
    \[\begin{tikzcd}
    	{\prod_{i\in I}\mathrm H\phi_i} & {\prod_{\set{i_0}}\mathrm H\phi_{i_0}} & \\
    	{\mathrm H\phi} & {\mathrm H(\set{i_0}\cdot\phi_{i_0})} & {\mathrm H\phi_{i_0}}
    	\arrow["{i_0^*}", from=1-1, to=1-2]
    	\arrow["{m_I}"', from=1-1, to=2-1]
    	\arrow["{m_{\{i_0\}}}"', from=1-2, to=2-2]
    	\arrow["\cong", from=1-2, to=2-3]
    	\arrow["{\mathrm H(i_0 \cdot \operatorname{id}_{\phi_i})}"', from=2-1, to=2-2]
    	\arrow["\cong"', from=2-2, to=2-3]
    \end{tikzcd}\]
    commutes. As the composite $\mathrm H\phi \to \mathrm H(\set{i_0} \cdot \phi_{i_0}) \to \mathrm H\phi_{i_0}$ coincides with the map $\mathrm H \iota_{i_0}$, it then follows that the composite 
    \[\begin{tikzcd}
    	{\prod_{i\in I}\mathrm H\phi_i} & {\mathrm H\phi} & {\mathrm H\phi_{i_0}}
    	\arrow["{m_I}", from=1-1, to=1-2]
    	\arrow["{\mathrm H \iota_{i_0}}", from=1-2, to=1-3]
    \end{tikzcd}\]
    is the canonical projection out of the product, and hence the composite 
    \[\begin{tikzcd}[sep = large]
    	{\prod_{i\in I}\mathrm H\phi_i} & {\mathrm H\phi} & {\prod_{i\in I}\mathrm H\phi_{i}}
    	\arrow["{m_I}", from=1-1, to=1-2]
    	\arrow["{(\mathrm H \iota_{i})_{i\in I}}", from=1-2, to=1-3]
    \end{tikzcd}\]
    is the identity.
    
    For the other composite, we consider the morphism $\mathsf p_I \colon I\cdot\phi\to\phi$ in $\FF/\sigma$ induced by the evident $\FF$-arrow $I \cdot \tau \to \tau$,\footnote{Note that $\mathsf p_I \colon I \cdot \phi \to \phi$ is the $I$-diagonal vf-arrow at $\phi$ in $(\FF/\sigma)^{\op}$ seen as a vf-category as in \Cref{ex:bbgamma-algebras}.} and we note that the composite
    \[\begin{tikzcd}
    	{\phi \cong I \cdot\phi_i} & {I \cdot \phi} & \phi
    	\arrow["{\operatorname{id} \cdot\, \iota_i}", from=1-1, to=1-2]
    	\arrow["{\mathsf p_I}", from=1-2, to=1-3]
    \end{tikzcd}\]
    is the identity of $\phi$ in $\FF/\sigma$.  By unitality and both naturality axioms, the diagram 
    \[\begin{tikzcd}[column sep = 45pt]
    	{\mathrm H\phi} & {\prod_{i\in I}\mathrm H\phi} & {\prod_{i\in I}\mathrm H\phi_{i}} \\
    	& {\mathrm H(I\cdot \phi)} & {\mathrm H(I\cdot \phi_i)}
    	\arrow["{\diag^{\mathrm H\phi}_I}", from=1-1, to=1-2]
    	\arrow["{\mathrm H\mathsf p_I}"', from=1-1, to=2-2]
    	\arrow["{\prod_{i\in I} \mathrm H \iota_i }", from=1-2, to=1-3]
    	\arrow["{m_I}", from=1-2, to=2-2]
    	\arrow["{m_I}", from=1-3, to=2-3]
    	\arrow["{\mathrm H(\operatorname{id}\cdot\, \iota_i)}"', from=2-2, to=2-3]
    \end{tikzcd}\]
    commutes, where $\diag^{\mathrm H\phi}_I$ is the $I$-diagonal at $\mathrm H{\phi}$ in $\Set$. Thus, the composite 
    \[\begin{tikzcd}[sep = large]
    	{\mathrm H\phi} & {\prod_{i\in I}\mathrm H\phi_{i}} & {\mathrm H(I\cdot \phi_i)}
    	\arrow["{(\mathrm H\iota_i)_{i\in I}}", from=1-1, to=1-2]
    	\arrow["{m_I}", from=1-2, to=1-3]
    \end{tikzcd}\]
    is also the identity. 
\end{proof}

Combining \Cref{prop:merg-struct} and \Cref{lem:merg-struct-preserve-coproducts}, we get the following. 

\begin{cor}[Merging inverts specialization]\label{lem:merging-is-inverse-to-specialization}
\looseness=-1
    Let $(\sigma_i)_{i\in I}$ be an $I$-family of filters. For each point $p$ and each $(I\cdot \sigma_i)$-family $(q_{i,s})_{(i,s): I\cdot\sigma_i}$ of points in $\X$, the map 
    \[\textstyle \X_{I \cdot \sigma_i} (p, q_{i,s} ) \to \prod_{i\in I} \X_{\sigma_i} (p, q_{i,s})\]
    defined by specialization along each coprojection $\sigma_i \hookrightarrow I\cdot \sigma_i$ is a bijection,    
    whose inverse sends an $I$-family $(u_i \colon p \vuto{\sigma_i} q_{i,s})_{i\in I}$ to their merge $\bigmrg_{i:I} u_i$.
\end{cor}

\begin{rmk}[Vf-arrows whose type is a set]\label{rmk:set-vf-arrows}
    It follows that, for any set $I$, the vf-arrows of type $I$ are fully determined by those of type $1$:
    \[\textstyle\X_I (p, q_i) \cong \prod_{i\in I} \X_1 (p, q_i).\]
    In particular, $\diag^p_0 \colon p\vuto{0}(\ )$ is the only vf-arrow of type $0$ with domain $p$.
\end{rmk}

\begin{rmk}
    In fact, \Cref{lem:merg-struct-preserve-coproducts} holds in far greater generality, as we sketch here. First note that, for a category $\mathcal C$, the full subcategory of $\bbgamma \mathcal C$ determined by filters of the form $\delta_S$ is the free completion of $\mathcal C$ under (small) products, $\operatorname{Prod}(\mathcal C)$: indeed, it is immediate to see that this subcategory coincides with $\operatorname{Fam}(\mathcal C^{\op})^{\op}$, where $\operatorname{Fam}$ is the well-known \emph{family construction} of \cite{benabouFiberedCategoriesFoundations1985}, freely adjoining (small) coproducts. Thus, there is an evident morphism of pseudomonads $\operatorname{Prod} \to \bbgamma$ which, adapting the proof of \cite[Lem.\ 3.1]{gambinoFormalTheoryPseudomonads2021} to lax morphisms, induces a forgetful 2-functor at the level of algebras in the opposite direction. In other words:
    \begin{enumerate}
        \item every $\bbgamma$-algebra is a $\operatorname{Prod}$-algebra, i.e.\ it admits small products;
        \item every lax morphism of $\bbgamma$-algebras is a lax morphism of $\operatorname{Prod}$-algebras, and by colax-idempotency every such morphism is actually a pseudomorphism (cf.\ \cite[Lem.\ 6.5 and Cor.\ 6.6]{kellyPropertylikeStructures1997}), so it preserves small products. 
    \end{enumerate}
\end{rmk}

\section{Tautness}\label{sec:tautness}

As we have seen in \Cref{ex:top-sp-as-vf-posets} and \Cref{rmk:relational-algebras}, the 2-category $\TopSp$ of topological spaces embeds both in vf-posets and in vu-posets, and the two embeddings commute with the forgetful 2-functor $\vfPos \to \vuPos$. 
\[\begin{tikzcd}[sep = small]
	{\mathsf{vfPos}} && \vuPos \\
	& \TopSp
	\arrow["{\mathsf U}", from=1-1, to=1-3]
	\arrow[hook', from=2-2, to=1-1]
	\arrow[hook, from=2-2, to=1-3]
\end{tikzcd}\]
As already noted in \cite[Rem.\ 3.14]{vangoolToposesEnoughPoints2026} and \cite[Thm.\ 5.5]{hamadGeneralisedUltracategoriesConceptual2025}, exploiting Barr's correspondence between topological spaces and relational algebras, the image of $\TopSp \mono \vuPos$ consists of those vu-posets $X$ satisfying the following condition:
\[ \mbox{if } p \vuleq{\tau} q_{\phi t} \mbox{ for some } \UF\mbox{-arrow } \phi \colon \tau \epi \sigma, \mbox{ then } p \vuleq{\sigma} q_{s}. \]
In \cite{saadiaGeometricMorphismsVirtual2026}, this property is studied under the name of \emph{tautness}, and it is extended to a condition on arbitrary vu-categories, which is satisfied in particular by the vu-category of points of any topos  (cf.\ \cite[Prop.\ 1.10]{saadiaGeometricMorphismsVirtual2026}).

\looseness=-1
In this section, we first recall the general notion of tautness for arbitrary vu-categories from \cite{saadiaGeometricMorphismsVirtual2026} and show how to adapt it to vf-categories. Next, continuing the above discussion of merging structures, we characterize tautness at the level of reindexing functors by introducing \emph{neighborhood structures}.
Exploiting this formalism, we prove that 
\begin{enumerate}
    \item tautness is preserved by the forgetful 2-functor $\vfCAT \to \vuCAT$ (\Cref{cor:U-preserves-tautness}), and that
    \item tautness precisely characterizes those vf-posets which arise from topological spaces (\Cref{cor:vf-poset=topsps}).
\end{enumerate}
In particular, this will entail that the forgetful 2-functor $\vfPos \to \vuPos$ induces an equivalence between taut vf-posets and taut vu-posets: this is the posetal case of the main result of this paper (\Cref{thm:main-theorem}), which will show that the same holds at the level of taut vf-\emph{categories} and taut vu-\emph{categories}.
We conclude this section by characterizing tautness as an orthogonality property (\Cref{prop:tautness-orthogonality-converse}).

\subsection{Definition and examples}\label{ssec:tautness-definition}
We begin by recalling the definition of tautness for vu-categories given in \cite{saadiaGeometricMorphismsVirtual2026}. 

\begin{defn}[Taut vu-category, \protect{\cite[Def.\ 1.6]{saadiaGeometricMorphismsVirtual2026}}]\label{def:taut-vu-cat}
    Let $\X$ be a vu-category. For a $\UF$-arrow $\phi \colon \tau \epi \sigma$, a vu-arrow $u\colon p\vuto{\tau}q_{\phi t}$ in $\X$ is \emph{reducible along $\phi$}, or \emph{$\phi$-reducible}, if $\psi_1^*u = \psi_2^*u$ holds for all pairs of $\UF$-arrows $\psi_1,\psi_2 \colon \lambda\rightrightarrows \tau$ such that $\phi \psi_1= \phi \psi_2$. We say that $\X$ is \emph{taut} if every $\phi$-reducible vu-arrow $u\colon p\vuto{\tau}q_{\phi t}$ can be uniquely reduced, in the sense that $u=\phi^*v$ for a unique vu-arrow $v\colon p\vuto{\sigma}q_s$.
\end{defn}

\begin{ex}
    Most vu-categories of interest are taut: the vu-category of points of a topos is always taut \cite[Prop.\ 1.10]{saadiaGeometricMorphismsVirtual2026}, and a vu-poset is taut precisely if it arises from a topological space \cite[Cor.\ 2.15]{saadiaGeometricMorphismsVirtual2026}.
\end{ex}

In order to extend the notion of tautness to vf-categories, we begin by observing that we cannot naïvely ask for arbitrary reducible vf-arrows to be uniquely reduced.

\begin{rmk}
    Let $\X$ be a vf-category and fix a point $p$ and a $\sigma$-family of points $(q_s)_{s:\sigma}$. For any two $u,u' \colon p \vuto{\sigma} q_s$, their reindexings along the canonical $\FF$-arrow $\phi \colon 0 \to \sigma$ coincide (see \Cref{rmk:set-vf-arrows}), and they are trivially $\phi$-reducible: therefore, if we were to require that every reducible vf-arrow can be uniquely reduced, it would follow that $u = u'$, i.e.\ that $\X_\sigma(p,q_s)$ is a singleton. 
\end{rmk}

The `correct analogue' of \Cref{def:taut-vu-cat} for vf-categories requires us to restrict the notion of reducibility to \emph{covers} in $\FF$ --- and, indeed, recall that all $\UF$-arrows are covers. Unlike in the case of vu-categories, the existence of kernel pairs in $\FF$ then allows us to simplify the definition as follows.

\begin{defn}[Taut vf-category]\label{defn:taut-vf-cat}
    Let $\X$ be a vf-category. For a cover $\phi \colon \tau \epi \sigma$ in $\FF$, a vf-arrow $u\colon p\vuto{\tau}q_{\phi t}$ in $\X$ is \emph{reducible along $\phi$}, or \emph{$\phi$-reducible}, if $\pi_1^*u = \pi_2^*u$ holds, where $\pi_1,\pi_2 \colon \tau \times_\sigma \tau \rightrightarrows \tau$ is the kernel pair of $\phi$. We say that $\X$ is \emph{taut} if every $\phi$-reducible vf-arrow $u\colon p\vuto{\tau}q_{\phi t}$ can be uniquely reduced, in the sense that $u=\phi^*v$ for a unique vf-arrow $v\colon p\vuto{\sigma}q_s$.
\end{defn}
\begin{rmk}\label{rmk:tautness-as-equalizer}
    In other words, a vf-category $\X$ is taut if for any point $p$, any $\sigma$-family of points $(q_s)_{s:\sigma}$, and any cover $\phi \colon \tau \epi \sigma$ in $\FF$, the diagram
\[\begin{tikzcd}
	{\X_\sigma(p,q_s)} & {\X_\tau(p,q_{\phi t})} & {\X_{\tau\times_\sigma\tau}(p,q_s)}
	\arrow["{\phi^*}", from=1-1, to=1-2]
	\arrow["{\pi_1^*}", shift left, from=1-2, to=1-3]
	\arrow["{\pi_2^*}"', shift right, from=1-2, to=1-3]
\end{tikzcd}\]
    is an equalizer in $\Set$, where $\pi_1,\pi_2 \colon \tau \times_\sigma \tau \rightrightarrows \tau$ is the kernel pair of $\phi$. 
\end{rmk}

\begin{rmk}[Reducibility and ultrafilters]\label{rmk:reducibility-ultrafilters}
    Let $u \colon p \vuto{\tau} q_{\phi t}$ be a vf-arrow in a vf-category $\X$, for a cover $\phi \colon \tau \epi \sigma$. A priori, $\phi$-reducibility is stronger than asking that $\psi_1^*u = \psi_2^*u$ for all pairs of $\FF$-arrows $\psi_1,\psi_2 \colon \lambda \rightrightarrows \tau$ such that $\phi\psi_1 = \phi\psi_2$ \emph{with $\lambda$ an ultrafilter}. However, we will show in \Cref{prop:vureduc-vfreduc} that the two notions of reducibility coincide if $\X$ is taut: this will be the main step towards proving that $\mathsf{U} \colon \vfCAT \to \vuCAT$ preserves tautness.
\end{rmk}

We shall now give some examples of taut vf-categories. First, exploiting the fact that a cover $\phi \colon \tau \epi \sigma$ in $\FF$ satisfies $\sigma = \phi_*\tau$ rather than merely $\sigma \subseteq \phi_*\tau$, the following can be shown exactly as in the case of vu-posets \cite[Prop.\ 2.6 and 2.7]{saadiaGeometricMorphismsVirtual2026}.

\begin{lem}\label{lem:top-spaces-are-taut}
    Topological spaces are taut as vf-posets.
\end{lem}

Similarly, the following can be proved exactly as its vu-categorical analogue \cite[Prop.\ 1.10]{saadiaGeometricMorphismsVirtual2026}.

\begin{prop}
    For any topos $\mathcal E$, the vf-category $\mathsf{pt}(\mathcal E)$ is taut.
\end{prop}

The first step in the proof of \cite[Prop.\ 1.10]{saadiaGeometricMorphismsVirtual2026} is to show that $\Set$ is a taut vu-category (cf.\ \cite[Lem.\ 1.9]{saadiaGeometricMorphismsVirtual2026}). We give here a direct proof that $\Set$ is taut as a vf-category, since the proof is much simpler than in the vu-categorical setting: indeed, it essentially follows from covers being effective epimorphisms in $\FF$ and from \Cref{cons:reduced-products}.

\begin{prop}\label{lem:set-is-taut}
    The vf-category $\Set$ is taut.
\end{prop}
\begin{proof}
    For a filter $\sigma$ on a set $S$, let $(B_s)_{s:\sigma}$ be a $\sigma$-family of sets and let $\phi \colon \tau \epi \sigma$ be a cover in $\FF$. Since $\phi$ is the coequalizer of its kernel pair $\pi_1,\pi_2 \colon \tau\times_\sigma\tau \rightrightarrows \tau$ (in $\FF$ and hence in $\FF/S$), the following diagram is an equalizer in $\Set$:
\[\begin{tikzcd}
	{\FF/S(\sigma, S\cdot B_s) } & {\FF/S(\tau, S\cdot B_s) } & {\FF/S(\tau\times_\sigma\tau, S\cdot B_s). }
	\arrow["{\phi^*}", from=1-1, to=1-2]
	\arrow["{\pi_1^*}", shift left, from=1-2, to=1-3]
	\arrow["{\pi_2^*}"', shift right, from=1-2, to=1-3]
\end{tikzcd}\]
    By \Cref{cons:reduced-products}, we can identify this diagram with the one below:
\[\begin{tikzcd}
	{\int_{s:\sigma} B_s} & {\int_{t:\tau} B_{\phi t}} & {\int_{\tau\times_\sigma\tau} B_s.}
	\arrow["{\phi^*}", from=1-1, to=1-2]
	\arrow["{\pi_1^*}", shift left, from=1-2, to=1-3]
	\arrow["{\pi_2^*}"', shift right, from=1-2, to=1-3]
\end{tikzcd}\]
    Thus, for any set $A$, the diagram below is an equalizer as well; recalling the description of vf-arrows in $\Set$ given in \Cref{ex:sets}, this concludes the proof by \Cref{rmk:tautness-as-equalizer}. 
    \[\begin{tikzcd}
	{\Set(A,\int_{s:\sigma} B_s)} & {\Set(A,\int_{t:\tau} B_{\phi t})} & {\Set(A,\int_{\tau\times_\sigma\tau} B_s).}
	\arrow["{\phi^*}", from=1-1, to=1-2]
	\arrow["{\pi_1^*}", shift left, from=1-2, to=1-3]
	\arrow["{\pi_2^*}"', shift right, from=1-2, to=1-3]
\end{tikzcd}\qedhere\]
\end{proof}

Consider now the $\bbgamma$-algebra $\mathsf{CompMet}_k$ of complete $k$-bounded metric spaces, for some $k \in (0,\infty]$, described in \Cref{ex:metric-spaces}; seeing it as a vf-category as in \Cref{ex:bbgamma-algebras}, we can follow the same proof strategy as for \Cref{lem:set-is-taut} to show the following.

\begin{prop}\label{lem:compmetk-is-taut}
    For $k \in (0,\infty]$, the vf-category $\mathsf{CompMet}_k$ is taut.
\end{prop}
\begin{proof}
    For a filter $\sigma$ on a set $S$, let $(B_s)_{s:\sigma}$ be a $\sigma$-family of complete $k$-bounded metric spaces, and let $\phi \colon \tau \epi \sigma$ be a cover in $\FF$. As in the proof of \Cref{lem:set-is-taut}, it suffices to show that the following diagram is an equalizer in $\mathsf{CompMet}_k$:
    \[\begin{tikzcd}
        {\int_{s:\sigma} B_s} & {\int_{t:\tau} B_{\phi t}} & {\int_{\tau\times_\sigma\tau} B_s.}
        \arrow["{\phi^*}", from=1-1, to=1-2]
        \arrow["{\pi_1^*}", shift left, from=1-2, to=1-3]
        \arrow["{\pi_2^*}"', shift right, from=1-2, to=1-3]
    \end{tikzcd}\]

    First, note that the category $\mathsf{CompMet}_k$ admits all equalizers; then, let $E$ be the equalizer of $\pi_1^*$ and $\pi_2^*$, which can be described as the (closed) subspace of $\int_{t:\tau}B_{\phi t}$ given by the points $(b_t)_{t:\tau}$ such that $\pi_1^* (b_t)_{t:\tau} = \pi_2^*(b_t)_{t:\tau}$.

    Since $\phi$ is a cover, the map $\phi^*$ is an \emph{isometric embedding}, i.e.\ an isometry onto its image. Indeed, for any $(x_s)_{s:\sigma}$ and $(y_s)_{s:\sigma}$ in $\int_{s:\sigma}B_s$, we can write:
    \begin{align*}
     d_\tau((x_{\phi t}),(y_{\phi t})) & = \inf_{T_0\in\tau}\sup_{t\in T_0}d_{B_{\phi t}}(x_{\phi t},y_{\phi t})\\ & = \inf_{T_0\in\tau}\sup_{s\in\phi[T_0]}d_{B_s}(x_s,y_s)\\ & = \inf_{S_0\in\sigma}\sup_{s\in S_0}d_{B_s}(x_s,y_s) = d_\sigma((x_s),(y_s))
    \end{align*}
    where the third equality follows from $\phi_*\tau = \sigma$.
    Moreover, it is clear that $\phi^*$ takes values in $E$, so that it induces an isometric embedding $\int_{s:\sigma}B_s\mono E$.
    
    To prove that the image of $\phi^*$ is precisely $E$, it suffices to show that its image is dense in $E$: indeed, the image of an isometric embedding between complete metric spaces is closed, and hence it coincides with the codomain if it is dense.
    Therefore, fix some $(z_t)_{t:\tau}$ in $E$ and some $\epsilon > 0$; we need to show that there exists some $(w_s)_{s:\sigma}$ in $\int_{s:\sigma}B_s$ such that $d_\tau ((z_t), (w_{\phi t})) \leq \epsilon$.
    
    The assumption $\pi_1^*(z_t) = \pi_2^*(z_t)$ gives $d_{B_{\phi t}}(z_t, z_{t'}) < \epsilon$ for $(t,t') : \tau\times_\sigma \tau$; as the kernel pair
    $\tau\times_\sigma\tau$ can be identified with the restriction of the product filter $\tau\times\tau$,
    this means that there exists some $T_0\in\tau$, contained in the domain of $(z_t)$, such that
    \[d_{B_{\phi t}}(z_t,z_{t'}) < \epsilon \quad \text{ for all } t,t'\in T_0 \text{ with } \phi t = \phi t'.\]
    Since $\phi$ is a cover, the set $S_0\coloneqq\phi[T_0]$ lies in $\sigma$: for each $s\in S_0$ pick some $t_s \in T_0$ such that $\phi(t_s) = s$, and let $(w_s)_{s:\sigma}$ be the $\sigma$-family of points determined by $(z_{t_s})_{s\in S_0}$. Since $d_{B_{\phi t}}(z_t, w_{\phi t}) < \epsilon$ for all $t \in T_0$, we conclude that $d_\tau ((z_t), (w_{\phi t})) \leq \epsilon$, as desired.
    The map $\int_{s:\sigma}B_s\to E$ is thus a bijective isometry, and hence an isomorphism in $\mathsf{CompMet}_k$.
\end{proof}

\begin{rmk}[What about $\mathsf{Met}_k$?]\label{rmk:metk-is-not-taut}
    Completeness is crucial in \Cref{lem:compmetk-is-taut} above; as we show now, the vf-category $\mathsf{Met}_k$ of $k$-bounded metric spaces is not taut. Fix a $k$-bounded metric space $X$ and take a proper filter $\sigma$, so that the unique $\phi\colon\sigma\epi\ptuf$ is a cover. Consider the diagram:
    \[\begin{tikzcd}
    	X & {X^\sigma} & {X^{\sigma\times \sigma}.}
    	\arrow["{\phi^*}", hook, from=1-1, to=1-2]
    	\arrow["{\pi_1^*}", shift left, from=1-2, to=1-3]
    	\arrow["{\pi_2^*}"', shift right, from=1-2, to=1-3]
    \end{tikzcd}\]

    On the one hand, the image of $X \hookrightarrow X^\sigma$ consists precisely of the $\sigma$-families of points $(x_s)_{s:\sigma}$ which converge to some point in $X$ for the topology induced by the metric. On the other hand, a $\sigma$-family $(x_s)_{s:\sigma}$ lies in the equalizer of $\pi_1^*$ and $\pi_2^*$ if and only if it is \emph{$\sigma$-Cauchy}, in the sense that for every $\epsilon >0$ there exists some $S_0 \in \sigma$ such that $d_X(x_s,x_{s'}) < \epsilon$ for all $s,s' \in S_0$. 
    For $\sigma$ the Fréchet filter on $\mathbb N$, it then follows that the above diagram is an equalizer if and only if every Cauchy sequence converges, i.e.\ if and only if $X$ is complete.
\end{rmk}

\begin{rmk}
    \looseness=-1
    Although it goes beyond the scope of this paper to prove it, we believe that $\mathsf{CompMet}_k$ is the \emph{tautification}, in a sense we will make precise in \Cref{rmk:tautification}, of $\mathsf{Met}_k$.  This suggests interpreting tautness as a \emph{local} (or \emph{internal}) Cauchy-completeness notion, saying that the points of the vf-category are themselves `complete'. 
\end{rmk}

We conclude this subsection by showing that tautness, for $\bbgamma$-algebras seen as vf-categories, cannot be axiomatized as an equational property. Indeed, on the one hand, we have just seen that there exist $\bbgamma$-algebras which are not taut --- namely, $\mathsf{Met}_k$; on the other hand, as we now show, free $\bbgamma$-algebras are always taut.

\begin{prop}
    For any category $\mathcal C$, the free $\bbgamma$-algebra $\bbgamma \mathcal C$ is taut as a vf-category.
    \begin{proof}
        Consider the multiplication $\bbmu_{\mathcal C} \colon \bbgamma^2 \mathcal C \to \bbgamma \mathcal C$ of the pseudomonad $\bbgamma$ at $\mathcal C$, giving the free $\bbgamma$-algebra structure; proceeding once again as in \Cref{lem:set-is-taut}, we need to show that, for any $\sigma$-family of filter-families $( (c_{s,t})_{t:\tau_s} )_{s:\sigma}$ in $\mathcal C$ and any cover $\phi \colon \lambda \epi \sigma$ in $\FF$, the diagram
\[\begin{tikzcd}
	{\bbmu_{\mathcal C}(( (c_{s,t})_{t:\tau_s} )_{s:\sigma}) } & {\bbmu_{\mathcal C}(( (c_{\phi \ell,t})_{t:\tau_{\phi \ell}} )_{\ell:\lambda})} & {\bbmu_{\mathcal C}(( (c_{s, t})_{t:\tau_{s}} )_{(\ell,s,\ell'):\lambda \times_\sigma \lambda})}
	\arrow["{\phi^*}", from=1-1, to=1-2]
	\arrow["{\pi_1^*}", shift left, from=1-2, to=1-3]
	\arrow["{\pi_2^*}"', shift right, from=1-2, to=1-3]
\end{tikzcd}\]
    is an equalizer in $\bbgamma \mathcal C$, where $\pi_1,\pi_2 \colon \lambda\times_\sigma\lambda \rightrightarrows \lambda$ is the kernel pair of $\phi$. Explicitly, the previous diagram can be rewritten as
\[\begin{tikzcd}
	{(c_{s,t})_{\sigma \cdot \tau_s}} & {(c_{\phi \ell,t})_{\lambda \cdot \tau_{\phi \ell}}} & {(c_{s,t})_{(\lambda\times_\sigma \lambda) \cdot \tau_s},}
	\arrow["{\phi^*}", from=1-1, to=1-2]
	\arrow["{\pi^*_1}", shift left, from=1-2, to=1-3]
	\arrow["{\pi^*_2}"', shift right, from=1-2, to=1-3]
\end{tikzcd}\]
    where the map $\phi^*$ essentially corresponds to the $\FF$-arrow $\phi \cdot \operatorname{id}_{\tau_s} \colon \lambda \cdot \tau_{\phi \ell} \to \sigma \cdot \tau_s$, and similarly $\pi_i^*$ corresponds to the $\FF$-arrow $\pi_i \cdot \operatorname{id}_{\tau_{\phi \ell}} \colon (\lambda \times_\sigma \lambda)\cdot \tau_s \to \lambda \cdot \tau_{\phi \ell}$. The fact that the above diagram is an equalizer can then be deduced from the fact, whose proof we leave to the reader, that the diagram below is a coequalizer in $\FF$.\footnote{More conceptually, this comes from the fact that the functor
    \[\FF/\sigma\to \FF/(\sigma\cdot\tau_s) \qquad (\phi\colon\lambda\to\sigma)\mapsto (\phi\cdot\operatorname{id}_{\tau_{\phi \ell}}\colon \lambda\cdot\tau_{\phi\ell} \to \sigma\cdot\tau_s)\]
    can be shown to have a right adjoint, corresponding to some \emph{dependent} internal-hom structure on $\FF$ which generalizes the closed structure of \Cref{rmk:monoidal-structure-on-F}; hence, it preserves colimits.}
\[\begin{tikzcd}[sep = 45pt]
	{(\lambda \times_\sigma \lambda)\cdot \tau_s} & {\lambda\cdot \tau_{\phi\ell}} & {\sigma\cdot \tau_s}
	\arrow["{\pi_1\cdot\operatorname{id}_{\tau_{\phi \ell}}}", shift left, from=1-1, to=1-2]
	\arrow["{\pi_2\cdot\operatorname{id}_{\tau_{\phi \ell}}}"', shift right, from=1-1, to=1-2]
	\arrow["{\phi\cdot \operatorname{id}_{\tau_s}}", from=1-2, to=1-3]
\end{tikzcd}\qedhere\]
    \end{proof}
\end{prop}

\subsection{Tautness as a sheaf condition}\label{ssec:tautness-sheaf}
As expressed in \cite{saadiaGeometricMorphismsVirtual2026} and recalled above, tautness can be conceived as a `tameness' condition on the behavior of reindexing.
In this subsection, we build on the merging structures of \Cref{ssec:merging} and characterize tautness of a vf-category in terms of its reindexing functors. 

First, recall from \Cref{ssec:category-of-filters} that, for each filter $\sigma$, the category $\FF/\sigma$ is regular and can thus be endowed with the regular topology, whose covering families are generated by covers. We can then observe that the condition of \Cref{rmk:tautness-as-equalizer} looks very similar to the sheaf condition described in \Cref{prop:sheaves-on-F-over-sigma}. This is not a coincidence: indeed, as we now show, tautness precisely corresponds to the reindexing functors being sheaves for the regular topology (cf.\ also \cite[Rem.\ 1.7]{saadiaGeometricMorphismsVirtual2026}). 

\begin{prop}[Tautness as a sheaf condition]\label{prop:tautness-as-sheaf-condition}
A vf-category $\X$ is taut if and only if, for each point $p$ and each $\sigma$-family of points $(q_s)_{s:\sigma}$, the reindexing functor
\[ (\FF/\sigma)^{\op} \to \Set \qquad (\phi \colon \tau \to \sigma) \mapsto \X_\tau(p, q_{\phi t}) \]
is a sheaf for the regular topology on $\FF/\sigma$.
\begin{proof}
    As we said, we will make use of the characterization of sheaves for the regular topology on $\FF/\sigma$ given in \Cref{prop:sheaves-on-F-over-sigma}. If each reindexing functor is a sheaf for the regular topology, then tautness follows from \Cref{rmk:tautness-as-equalizer} by noting that every cover $\phi \colon \tau \epi \sigma$ in $\FF$ induces a cover $\phi \epi \operatorname{id}_\sigma$ in $\FF/\sigma$.
    
    Conversely, suppose that $\X$ is taut. We need to show that, for any cover $\theta \colon \psi \epi \phi$ in $\FF/\sigma$, the diagram 
    \[\begin{tikzcd}
        {\X_\lambda(p, q_{\phi \ell})} & {\X_{\kappa}(p, q_{\psi k})} & {\X_{\kappa\times_{\lambda}\kappa}(p, q_{\ell})}
        \arrow["{\theta^*}", from=1-1, to=1-2]
        \arrow["{\pi_1^*}", shift left, from=1-2, to=1-3]
        \arrow["{\pi_2^*}"', shift right, from=1-2, to=1-3]
    \end{tikzcd}\]
    is an equalizer in $\Set$, where the diagram below represents the kernel pair of $\theta$ in $\FF/\sigma$.
    \[\begin{tikzcd}
        {\kappa\times_{\lambda}\kappa} & {\kappa} & \lambda \\
        & \sigma
        \arrow["{\pi_1}", shift left, from=1-1, to=1-2]
        \arrow["{\pi_2}"', shift right, from=1-1, to=1-2]
        \arrow["{\psi\times_\phi\psi}"', curve={height=12pt}, from=1-1, to=2-2]
        \arrow["\theta", two heads, from=1-2, to=1-3]
        \arrow["\psi", from=1-2, to=2-2]
        \arrow["\phi", curve={height=-12pt}, from=1-3, to=2-2]
    \end{tikzcd}\]
    We take some vf-arrow $u \colon p \vuto{\kappa} q_{\psi k}$ such that $\pi_1^* u = \pi_2^* u$. Since $\pi_1,\pi_2$ is the kernel pair of $\theta$ in $\FF$, the vf-arrow $u \colon p \vuto{\kappa}q_{\phi(\theta k)}$ is $\theta$-reducible. Thus, by tautness, there exists a unique vf-arrow $v \colon p \vuto{\lambda} q_{\phi \ell}$ in $\X$ such that $u = \theta^* v$, which concludes the proof.
    \end{proof}
\end{prop}

Expanding our discussion on merging in \Cref{ssec:merging}, we can thus give the following definition, which allows us to rephrase \Cref{prop:tautness-as-sheaf-condition} in the language of merging structures.

\begin{defn}[Neighborhood structures]\label{defn:neighborhood-structure}
    A \emph{neighborhood structure} is a merging structure $\mathrm H \colon (\FF/\sigma)^{\op}\to\Set$ 
    which is also a sheaf for the regular topology on $\FF/\sigma$.
\end{defn}

\begin{cor}\label{cor:taut-neighbourhood}
    A vf-category is taut if and only if its reindexing functors, seen as merging structures, are in fact neighborhood structures.
\end{cor}

\begin{rmk}
    Recall from \Cref{rmk:F-is-coherent} that we can also consider the coherent topology on $\FF/\sigma$ and, as shown in \Cref{lem:merg-struct-preserve-coproducts}, reindexing functors map finite coproducts to products. Hence, a merging structure is a neighborhood structure also if and only if it is a sheaf for the coherent topology, and similarly, a vf-category is taut if and only if its reindexing functors are sheaves for the coherent topology.
    
    As the reindexing functors map \emph{all} (small) coproducts to products, one might be tempted to extend the above remark to the infinitary topology $J_{\infty}$ considered in \cite[\S 3]{awodeyUltrasheavesDoubleNegation2004}. This is, however, not possible as the latter is not subcanonical, contrary to what is stated in \cite[Lem.\ 3.3]{awodeyUltrasheavesDoubleNegation2004}, essentially since $\FF$ is not infinitary extensive (cf.\ \Cref{prop:F-infextensive}).
\end{rmk}

\subsection{Tautness from vf-categories to vu-categories}\label{ssec:tautness-preserved}
\looseness=-1
In this subsection, we show that tautness is preserved by the forgetful 2-functor $\mathsf U\colon \vfCAT \to \vuCAT$. The main ingredient of the proof can be stated, in terms of neighborhood structures, as the following lemma.

\begin{lem}\label{lem:neighbourhood-ultra}
    Let $\mathrm H \colon (\FF/\sigma)^{\op} \to \Set$ be a neighborhood structure and let $\phi \colon \tau \to \sigma$ be an $\FF$-arrow.    For all $u,v \in \mathrm H\phi$, if $\mathrm H \iota (u) = \mathrm H \iota (v)$ holds for each ultrafilter $\tau'$ refining $\tau$, with the evident $\FF$-arrow $\iota \colon \tau' \hookrightarrow \tau$, then $u = v$. 
\end{lem}
\begin{proof}
    Let $L$ be a set indexing the ultrafilters $(\iota_\ell \colon \tau_\ell\hookrightarrow \tau)_{\ell \in L}$ refining $\tau$. Recall from \Cref{lem:cover-by-uf} that the canonical $\FF$-arrow $\psi \colon L\cdot\tau_\ell \epi \tau$ is a cover, which induces a morphism $L \cdot \phi_\ell \epi \phi$ in $\FF/\sigma$, where $\phi_\ell \coloneqq \phi \circ \iota_\ell \colon \tau_\ell \to \sigma$. As in the proof of \Cref{lem:merg-struct-preserve-coproducts}, the diagram 
    \[\begin{tikzcd}[column sep = 45 pt]
    	{\mathrm H\phi} & {\prod_{\ell \in L}\mathrm H\phi} & {\prod_{\ell \in L}\mathrm H\phi_{\ell}} \\
    	& {\mathrm H(L\cdot \phi)} & {\mathrm H(L\cdot \phi_\ell)}
    	\arrow["{\diag^{\mathrm H\phi}_L}", from=1-1, to=1-2]
    	\arrow["{\mathrm H \mathsf p_L}"', from=1-1, to=2-2]
    	\arrow["{\prod_{\ell\in L} \mathrm H \iota_\ell}", from=1-2, to=1-3]
    	\arrow["{m_L}", from=1-2, to=2-2]
    	\arrow["{m_L}", from=1-3, to=2-3]
    	\arrow["{\mathrm H(\operatorname{id}\cdot\, \iota_\ell)}"', from=2-2, to=2-3]
    \end{tikzcd}\]
    commutes, where the morphism $\mathsf p_L \colon L \cdot \phi \to \phi$ is induced by the evident $\FF$-arrow $L\cdot \tau \to \tau$ and $\diag^{\mathrm H\phi}_L \colon \mathrm H\phi \to \prod_{\ell\in L}\mathrm H\phi$ is the diagonal map. Since the bottom composite coincides with the map $\mathrm H \psi $ and the one on top coincides with the canonical map $(\mathrm H \iota_\ell)_{\ell\in L}$, we can rewrite this diagram as below.
    \[\begin{tikzcd}[column sep = 45pt]
    	{\mathrm H\phi} & & {\prod_{\ell \in L}\mathrm H\phi_\ell} \\
    	& & {\mathrm H(L \cdot \phi_\ell)}
    	\arrow["{(\mathrm H \iota_\ell )_{\ell\in L}}", from=1-1, to=1-3]
    	\arrow["{\mathrm H \psi}"', from=1-1, to=2-3]
    	\arrow["{m_L}", from=1-3, to=2-3]
    \end{tikzcd}\]
    To conclude, note that the map $m_L$ is bijective by \Cref{lem:merg-struct-preserve-coproducts}, while the map $\mathrm H\psi$ is injective by the sheaf condition, as $\psi$ is a cover, so that the top map is injective as well.
\end{proof}

\begin{cor}\label{cor:vureduc-ultrafilters}
    Let $\X$ be a taut vf-category and let $u,v \colon p \vuto{\sigma} q_s$ be parallel vf-arrows. If $\iota^* u = \iota^* v$ for each ultrafilter $\tau$ refining $\sigma$ with evident $\FF$-arrow $\iota \colon \tau \hookrightarrow \sigma$, then $u = v$.
\begin{proof}
    Combine \Cref{lem:neighbourhood-ultra} and \Cref{cor:taut-neighbourhood}. 
\end{proof}
\end{cor}

\begin{prop}[Ultrafilters detect reducibility]\label{prop:vureduc-vfreduc}
    Let $\X$ be a taut vf-category and let $\phi \colon \tau \epi \sigma$ be a cover in $\FF$. A vf-arrow $u \colon p \vuto{\tau} q_{\phi t}$ is $\phi$-reducible if and only if $\psi_1^* u = \psi_2^* u$ for all pairs of $\FF$-arrows $\psi_1,\psi_2 \colon \lambda \rightrightarrows \tau$ such that $\phi\psi_1 = \phi\psi_2$ with $\lambda$ an ultrafilter.
\end{prop}
\begin{proof}
    The condition is clearly necessary; we show it is also sufficient. Let $\pi_1,\pi_2 \colon \tau\times_\sigma\tau \rightrightarrows \tau$ be the kernel pair of $\phi$.
    For any ultrafilter $\iota \colon \lambda \mono \tau\times_\sigma\tau$ refining $\tau\times_\sigma\tau$, we have
    \[\iota^*\pi_1^*u = (\pi_1\iota)^*u = (\pi_2\iota)^*u  = \iota^*\pi_2^*u,\]
    from which we conclude by \Cref{cor:vureduc-ultrafilters}.
\end{proof}

As hinted at in \Cref{rmk:reducibility-ultrafilters}, the previous proposition entails that, in a taut vf-category $\X$, reducibility of a vf-arrow along a $\UF$-arrow is equivalent to its reducibility as a vu-arrow in the vu-category $\mathsf{U}\X$. This immediately yields the following corollary.

\begin{cor}[$\mathsf U$ preserves tautness]\label{cor:U-preserves-tautness}
    If $\X$ is a taut vf-category, then $\mathsf U\X$ is a taut vu-category.
\end{cor}
\begin{proof}

    Let $\X$ be a taut vf-category, let $\phi\colon \tau\coverto \sigma$ be a $\UF$-arrow and let $u\colon p\vuto{\tau}q_{\phi t}$ be a vf-arrow in $\X$ that is $\phi$-reducible as a vu-arrow in $\mathsf U\X$, meaning that $\psi_1^*u = \psi_2^*u$ for all pairs of $\UF$-arrows $\psi_1,\psi_2 \colon \lambda \rightrightarrows \tau$ such that $\phi\psi_1 = \phi\psi_2$. \Cref{prop:vureduc-vfreduc} then ensures that $u$ is $\phi$-reducible as a vf-arrow in $\X$ and hence, by tautness of $\X$, there is a unique $v\colon p\vuto{\sigma}q_s$ such that $\phi^*v = u$. 
\end{proof}

\subsection{Taut vf-posets are topological spaces}\label{ssec:taut-vf-posets}Proceeding as in \cite[\S 2]{saadiaGeometricMorphismsVirtual2026} for the case of vu-posets, in this subsection we show that the notion of tautness precisely characterizes those vf-posets arising from topological spaces. To do this, we will consider the \emph{upset topology} on a vf-poset, first introduced in \cite[Def.\ 3.23]{vangoolToposesEnoughPoints2026} for vu-posets, whose opens are the `upsets' with respect to the convergence relation. 

\begin{defn}[Upset topology]\label{defn:upset-topology}
    The \emph{upset topology} of a vf-poset $X$ is the topology on the set $X$ whose opens, called \emph{upsets}, are subsets $U$ such that, for each vf-inequality $p\vuleq{\sigma} q_s$ in $X$,
    \[p\in U \Rightarrow \forall_{s:\sigma}\ q_s \in U.\]
\end{defn}

The proof of the following is straightforward (cf.\ \cite[Prop.\ 2.19]{saadiaGeometricMorphismsVirtual2026}).

\begin{prop}\label{prop:topref-adj}
    The upset topology defines a 2-functor $\Tref{-} \colon \vfPos\to\TopSp$ which is left 2-adjoint to the embedding $\TopSp\hookrightarrow\vfPos$.
\end{prop}

We will show that this adjunction restricts to an equivalence $\TopSp\simeq\vfPos_{\mathrm{taut}}$ between topological spaces and taut vf-posets. The key ingredient of the proof is given by \emph{neighborhood filters}: recall that for a point $p$ of a topological space $X$, the set of neighborhoods of $p$ forms a filter on $X$, which we denote by $\lambda_p$. This filter is characterized by the property that, for any $\sigma$-family of points $(q_s)_{s:\sigma}$ in $X$, we have $p \vuleq{\sigma} q_s$ if and only if $q_*\sigma$ refines $\lambda_p$. Below, we construct this filter in the language of taut vf-posets and neighborhood structures, and we employ it to derive the desired equivalence.

As we are working with vf-posets rather than arbitrary vf-categories, we restrict our attention to merging structures $\mathrm H \colon (\FF/\sigma)^{\op}\to \Set$ which are \emph{subterminal}, i.e.\ such that the set $\mathrm H(\phi)$ has at most one element for each $\FF$-arrow $\phi \colon \tau \to \sigma$.

\begin{prop}[Subterminal neighborhood structures are representable]\label{prop:neighborhood-porp}
    Let $\mathrm H \colon (\FF/\sigma)^{\op}\to \Set$ be a subterminal neighborhood structure. Then, there exists a refinement $\lambda$ of $\sigma$ such that $\mathrm H$ is represented by the monomorphism $\lambda\mono \sigma$.
\end{prop}
\begin{proof}

Let $L$ be a set indexing the refinements $(\iota_\ell \colon \sigma_\ell \hookrightarrow \sigma)_{\ell \in L}$ of $\sigma$ such that $\mathrm H \iota_\ell$ is nonempty. Consider the canonical $\FF$-arrow $\psi \colon L \cdot \sigma_\ell \to \sigma$ and, by regularity of $\FF$, let $\iota \colon \lambda \hookrightarrow \sigma$ be its image. To prove that $\iota$ represents $\mathrm H$ we need to show that, for each $\FF$-arrow $\phi \colon \tau \to \sigma$, the set $\mathrm H \phi$ is nonempty if and only if $\phi$ factors (uniquely) through $\iota$.
    \[\begin{tikzcd}
    	\tau & \\
    	\lambda & \sigma
    	\arrow[dashed, from=1-1, to=2-1]
    	\arrow["\phi", curve={height=-12pt}, from=1-1, to=2-2]
    	\arrow["\iota"', hook, from=2-1, to=2-2]
    \end{tikzcd}\]
    
First, suppose that such a factorization exists. By \Cref{lem:merg-struct-preserve-coproducts}, we have $\mathrm H\psi \cong \prod_{\ell\in L}\mathrm H\iota_\ell$ and so $\mathrm H\psi$ is nonempty; moreover, considering the cover $\psi \epi \iota$ in $\FF/\sigma$, the sheaf condition ensures that $\mathrm H\iota$ is nonempty as well, and hence so is $\mathrm H\phi$ by the factorization.

Conversely, suppose that $\mathrm H \phi$ is nonempty. Consider the image $\iota' \colon \phi_*\tau \hookrightarrow \sigma$ of $\phi$; by the sheaf condition, $\mathrm H\iota' $ is nonempty, so that $\iota ' = \iota_{\ell_0}$ for some ${\ell_0} \in L$. The coprojection $\sigma_{\ell_0} \mono L\cdot \sigma_\ell$ induces the desired factorization of $\phi$ through $\iota$ as follows.
    \[\begin{tikzcd}[baseline=(current bounding box.south)]
    	\tau && \\
    	{\sigma_{\ell_0}} \\
    	{L\cdot\sigma_\ell} & \lambda & \sigma
    	\arrow[two heads, from=1-1, to=2-1]
    	\arrow["\phi", curve={height=-12pt}, from=1-1, to=3-3]
    	\arrow[dashed, hook, from=2-1, to=3-1]
    	\arrow["\iota_{\ell_0}"{description}, curve={height=-6pt}, hook, from=2-1, to=3-3]
    	\arrow[two heads, from=3-1, to=3-2]
    	\arrow["\iota"', hook, from=3-2, to=3-3]
    \end{tikzcd}\qedhere\]
\end{proof}

\begin{cons}\label{cor:lambda_p}
    Let $X$ be a taut vf-poset and let $p$ be a point of $X$. Applying \Cref{prop:neighborhood-porp} to the reindexing functor at $p$ and the $\delta_X$-family $(x)_{x\in X}$ defined by the identity map, that is, $(q \colon \sigma\to X) \mapsto p\vuleq{\sigma} q_s$,  
    we obtain a filter $\lambda_p$ on the set $X$ such that a vf-inequality $p\vuleq{\sigma}q_s$ holds if and only if $q_*\sigma$ refines $\lambda_p$. 
\end{cons}

Thus, simply by reasoning on the neighborhood structure of the reindexing functors of a taut vf-poset $X$, we have constructed, for each point $p$ of $X$, a filter $\lambda_p$ which has the characterizing property of the neighborhood filter. However, to show that it indeed coincides with the filter of neighborhoods in the upset topology, we will need to make explicit use of the vf-categorical structure of $X$.

\begin{prop}[Filter of neighborhoods]\label{prop:lambda_p}
    Let $X$ be a taut vf-poset. For each point $p$, the filter $\lambda_p$ coincides with the filter of neighborhoods of $p$ in the upset topology.
\end{prop}
\begin{proof}
    Let $A$ be a neighborhood of $p$ in the upset topology and let $U$ be an upset containing $p$ and contained in $A$. Since $\lambda_p$ refines itself, we have, by \Cref{cor:lambda_p}, that $p \vuleq{\lambda_p} x$, where $(x)_{x:\lambda_p}$ is the $\lambda_p$-family defined by the identity map $\lambda_p \hookrightarrow X$. It follows that $U \in \lambda_p$ and hence $A \in \lambda_p$ as well.     
    Conversely, let $A \in \lambda_p$ and consider the subset $U\subseteq X$ defined by
    \[ U \coloneqq \set{ x \in X | A \in \lambda_x}.\]
    It is clear that $p$ lies in $U$; we will show that $U\subseteq A$, and that $U$ is an upset.

    Let $x\in U$, meaning that $A\in\lambda_x$. By \Cref{cor:lambda_p}, the vf-inequality $x \vuleq{1} x$ entails that $x_*1$ refines $\lambda_x$, where $x_*1$ is the principal ultrafilter at $x$; hence, $A$ is $x_*1$-large, meaning that $x \in A$.
    To show that $U$ is an upset, suppose that $a \vuleq{\sigma} b_s$ with $a \in U$ and consider the following commutative diagram in $\FF$.
    \[\begin{tikzcd}[sep = 16pt]
    	& {\sigma\cdot A} && A \\
    	{\sigma\cdot\lambda_{b_s}} && {\lambda_a} \\
    	& {\sigma\cdot X} && X
    	\arrow[from=1-2, to=1-4]
    	\arrow[hook, from=1-2, to=3-2]
    	\arrow["\lrcorner"{anchor=center, pos=0.125}, draw=none, from=1-2, to=3-4]
    	\arrow[hook, from=1-4, to=3-4]
    	\arrow[dotted, from=2-1, to=1-2]
    	\arrow[dashed, from=2-1, to=2-3]
    	\arrow[from=2-1, to=3-2]
    	\arrow[dashed, from=2-3, to=1-4]
    	\arrow[hook, from=2-3, to=3-4]
    	\arrow[from=3-2, to=3-4]
    \end{tikzcd}\]
    The dashed arrow $\lambda_a \to A$ is given by the fact that $A \in \lambda_a$ since $a \in U$, while the one $\sigma \cdot \lambda_{b_s} \to \lambda_{a}$ is witnessed by \Cref{cor:lambda_p}, using the composite $a\vuleq{\sigma}(b_s\vuleq{\lambda_{b_s}} x)$. Finally, the dotted arrow $\sigma\cdot \lambda_{b_s} \to \sigma\cdot A$ is induced by the universal property of the pullback, and it implies that $A \in \lambda_{b_s}$ for $s : \sigma$; hence, $b_s \in U$ for $s:\sigma$.
\end{proof}

\begin{cor}[Convergence in the upset topology]\label{cor:vf-cvg=upset-cvg}
    For a taut vf-poset $X$, a vf-inequality $p\vuleq{\sigma} q_s$ holds in $X$ if and only if it holds in $\Tref{X}$.
\end{cor}
\begin{proof}
    By \Cref{cor:lambda_p}, a vf-inequality $p\vuleq{\sigma} q_s$ holds in $X$ if and only if $q_*\sigma$ refines $\lambda_p$, which by \Cref{prop:lambda_p} is equivalent to $q_*\sigma$ converging to $p$ in the upset topology, i.e.\ to $p\vuleq{\sigma} q_s$ holding in $\Tref{X}$.
\end{proof}

\begin{cor}[Taut vf-posets are topological spaces]\label{cor:vf-poset=topsps}
    The embedding $\TopSp \hookrightarrow \vfPos$ corestricts to a 2-equivalence between topological spaces and taut vf-posets.
\end{cor}
\begin{proof}
    Immediate from \Cref{cor:vf-cvg=upset-cvg}, recalling from \Cref{lem:top-spaces-are-taut} that topological spaces are taut as vf-posets. 
\end{proof}

\subsection{Tautness as an orthogonality property} \label{ssec:orthogonality}
In the previous subsection, we have seen that --- exactly as in the vu-categorical case --- taut vf-posets can be identified with topological spaces. This justifies our notion of tautness at the posetal level. In this subsection, we characterize tautness for vf-categories as an \emph{orthogonality} property (\Cref{prop:tautness-orthogonality-converse}), induced by tautness for vf-posets (and similarly for vu-categories and vu-posets): in some sense, this allows us to understand \Cref{defn:taut-vf-cat} as the minimal condition generated by posetal tautness.

To introduce the maps defining the orthogonality condition, we begin with the following observation. Through the identification of topological spaces with taut vf-posets, we can see the construction of the upset topology as a \emph{tautification} operation on vf-posets, $\Tref{-} \colon \vfPos \to \vfPos_{\mathrm{taut}}$. In the following lemma, we describe the tautification $\Tref X$ of a vf-poset $X$ in terms of its vf-inequalities.

\begin{lem}[Tautification of vf-posets]\label{lem:tautification}
    Let $X$ be a vf-poset. A vf-inequality $p\vuleq{\sigma}q_s$ holds in $\Tref{X}$ if and only if there exists a cover $\phi \colon \tau \coverto \sigma$ such that $p\vuleq{\tau}q_{\phi t}$ holds in $X$.
\begin{proof}
    Let $X'$ be the vf-poset with the same points as $X$ and where a vf-inequality $p\vuleq{\sigma}q_s$ holds if and only if there exists a cover $\phi \colon \tau \coverto \sigma$ such that $p\vuleq{\tau}q_{\phi t}$ holds in $X$. It is straightforward to see that $X'$ is taut, hence a topological space, and that its topology coincides with the upset topology of $X$, so that $X' = \Tref{X}$.
\end{proof}
\end{lem}

\begin{defn}[Tautification maps]\label{defn:tautification}
    For each vf-poset $X$, the unit of the adjunction between $\Tref{-}$ and the inclusion $\vfPos_{\mathrm{taut}} \hookrightarrow \vfPos$ at $X$ is a vf-functor 
    \[ \rho_X \colon X \to \Tref{X}\]
    acting as the identity both on points and on vf-arrows: we refer to $\rho_X$ as the \emph{tautification map} associated to $X$. 
\end{defn}

\begin{prop}\label{prop:tautness-orthogonality-vf}
    Let $\X$ be a taut vf-category. Then, for any vf-poset $Y$, $\X$ is right orthogonal to the tautification map $\rho_Y \colon Y\to\Tref{Y}$, i.e.\ the induced map
    \[\rho_Y^* \colon \vfCAT(\Tref{Y},\X)\to\vfCAT(Y,\X)\]
    is an isomorphism of categories.
\end{prop}
\begin{proof}
    We first show that $\rho_Y^*$ is fully faithful.
    Faithfulness of $\rho_Y^*$ follows from $\rho_Y$ being the identity on points. 
    Now let $F , G \colon \Tref{Y} \to \X$ be vf-functors and let $\alpha \colon F\rho_Y \Rightarrow G\rho_Y$ be a vf-transformation whose components are given by $\alpha_y \colon F(y) \vuto{\ptuf} G(y)$ for each $y \in Y$; we need to show that $\alpha$ lifts to a vf-transformation $F\Rightarrow G$, i.e.\ that
    \[ Gu \circ_\ptuf \alpha_p = \alpha_{q_s} \circ_\sigma Fu\]
    holds for each vf-inequality $u \colon p \vuleq{\sigma} q_s$ in $\Tref{Y}$.
    By \Cref{lem:tautification}, $u$ is witnessed by a vf-inequality $v \colon p \vuleq{\tau} q_{\phi t}$ in $Y$ for some cover $\phi \colon \tau \epi \sigma$ with $\phi^*u = \rho_Y(v)$.
    Hence, using the naturality of $\alpha : F\rho_Y \Rightarrow G\rho_Y$, we have
    \begin{multline*}
     \phi^* (Gu \circ_\ptuf \alpha_p ) = \phi^* (Gu) \circ_\ptuf \alpha_p = G(\phi^*u) \circ_\ptuf \alpha_p = G\rho_Y(v) \circ_\ptuf \alpha_p = \\ = \alpha_{q_{\phi t}} \circ_\tau F\rho_Y(v) =\alpha_{q_{\phi t}}\circ_\tau F(\phi^* u) = \alpha_{q_{\phi t}} \circ_\tau \phi^* (Fu) = \phi^* ( \alpha_{q_s}\circ_\sigma Fu),
    \end{multline*}
    whence, by tautness, $Gu \circ_\ptuf \alpha_p = \alpha_{q_s} \circ_\sigma Fu$.

    We then show that $\rho_Y^*$ is bijective on objects: fix a vf-functor $F \colon Y \to \X$. We define an extension $\widetilde F \colon \Tref{Y} \to \X$ as follows.
    On points, we set $\widetilde F(y) \coloneqq F(y)$ for each $y \in Y$.
    Now let $u\colon p \vuleq{\sigma} q_s$ be a vf-inequality in $\Tref{Y}$, again witnessed by a vf-inequality $v \colon p \vuleq{\tau} q_{\phi t}$ in $Y$ for some cover $\phi \colon \tau \epi \sigma$. By posetality of $Y$ and by vf-functoriality of $F$, the vf-arrow $F(v) \colon F(p) \vuto{\tau} F(q_{\phi t})$ is $\phi$-reducible; therefore, by tautness of $\X$, there exists a unique vf-arrow $w \colon F(p) \vuto{\sigma} F(q_s)$ such that $\phi^* w = F(v)$, so that we can set $\widetilde F(u) \coloneqq w$.
    It is straightforward to check that this $w$ does not depend on the choice of $v$, that this construction defines a vf-functor $\widetilde F \colon \Tref{Y} \to \X$, and that this vf-functor is the unique one such that $\widetilde F \rho_Y = F$.
\end{proof}

\begin{rmk}[Tautification]\label{rmk:tautification}
    The previous proposition shows that the association $Y \mapsto \Tref{Y}$ defines a left 2-adjoint to the inclusion $\vfCAT_{\mathrm{taut}} \hookrightarrow \vfCAT$ \emph{relative} to the inclusion of vf-posets, in the sense of, e.g., \cite[Def.\ 2.9]{altenkirchMonadsNeedNot2015}.
\[\begin{tikzcd}[sep = small]
	& {\vfCAT_{\mathrm{taut}}} & \\
	\vfPos && \vfCAT
	\arrow[""{name=0, anchor=center, inner sep=0}, hook, from=1-2, to=2-3]
	\arrow[""{name=1, anchor=center, inner sep=0}, "{\Tref{-}}", from=2-1, to=1-2]
	\arrow[hook, from=2-1, to=2-3]
	\arrow["\dashv"{anchor=center}, draw=none, from=1, to=0]
\end{tikzcd}\]
\looseness=-1
We expect that a left 2-adjoint to the inclusion $\vfCAT_{\mathrm{taut}} \hookrightarrow \vfCAT$ exists \emph{tout court}, specializing to $Y \mapsto \Tref Y$ for vf-posets, but we leave this for future research.
\end{rmk}

\looseness=-1
With a similar definition of the tautification for vu-posets (cf.\ \cite[\S 2]{saadiaGeometricMorphismsVirtual2026}), and the same argument as in the above proof, one can also show a vu-categorical analogue of \Cref{prop:tautness-orthogonality-vf}.

\begin{prop}\label{prop:tautness-orthogonality}
    Taut vu-categories are right orthogonal to all tautification maps $\rho_Y \colon Y\to\Tref{Y}$ with $Y$ a vu-poset.
\end{prop}

In fact, the converse of \Cref{prop:tautness-orthogonality-vf} (and similarly of \Cref{prop:tautness-orthogonality}) is also true, hence characterizing tautness as an orthogonality property. For narrative purposes, we prove this fact here, even though the proof uses notions that will only be introduced in \Cref{sec:equivalence}; this is not circular, however, as only \Cref{prop:tautness-orthogonality} will be used in the proof of the main theorem.

\begin{prop}[Tautness as an orthogonality property]\label{prop:tautness-orthogonality-converse}
    Let $\X$ be a vf-category. Then, $\X$ is taut if and only if $\X$ is right orthogonal to all tautification maps $\rho_Y \colon Y\to\Tref{Y}$ with $Y$ a vf-poset.
\end{prop}
\begin{proof}
    One direction is given by \Cref{prop:tautness-orthogonality-vf}, so suppose that $\X$ is right orthogonal to all tautification maps $\rho_Y \colon Y\to\Tref{Y}$ with $Y$ vf-poset. To show that $\X$ is taut, fix a cover $\phi\colon \tau \epi \sigma$, represented by some map $\phi \colon T \epi S$, and consider the vf-poset $\mathrm I_{S,\phi}$ of \Cref{cons:walking-phi-reduc}.
    
    Fix points $p$ and $(q_s)_{s\in S}$ in $\X$.
    As explained in \Cref{rmk:walking-phi-reduc}, a vf-functor $F\colon \mathrm I_{S,\phi}\to\X$ mapping $a\mapsto p$ and $b_s\mapsto q_s$ for each $s\in S$ corresponds to a $\phi$-reducible vf-arrow $u \colon p\vuto{t:\tau}q_{\phi t}$.
    Moreover, $\Tref{\mathrm I_{S,\phi}}$ recovers the walking vf-arrow $\mathrm I_{S,\sigma}$ of \Cref{cons:walking-vf-arrow}; hence, an extension $\widetilde F \colon \Tref{\mathrm I_{S,\phi}}\to \X$ of such an $F$ uniquely corresponds to a vf-arrow $v\colon p\vuto{s:\sigma}q_s$ satisfying $\phi^*v = u$. In other words, the orthogonality property along $\rho_{\mathrm I_{S,\phi}} \colon \mathrm I_{S,\phi} \to \Tref{\mathrm I_{S,\phi}}$ precisely means that any $\phi$-reducible vf-arrow can be uniquely reduced.
\end{proof}
\begin{rmk}
    Clearly, the same argument also shows a converse to \Cref{prop:tautness-orthogonality} in the vu-categorical setting.
\end{rmk}

\section{The equivalence}\label{sec:equivalence}

In this last section, we finally prove the main result of this paper: the forgetful 2-functor $\mathsf U\colon \vfCAT \to \vuCAT$ restricts to a 2-equivalence of 2-categories between taut vf-categories and taut vu-categories. As we discussed in the introduction, this result can be understood as a categorification of the observation that topological convergence relations of filters are uniquely determined by their restrictions to ultrafilters.

Our strategy in this section is based on introducing \emph{walking vf-arrows} (\Cref{cons:walking-vf-arrow}), a collection of topological spaces suitably classifying vf-arrows of each type. This will allow us to endow the class of points of a taut vu-category with a family of vf-arrows extending its original vu-arrows (\Cref{thm:obtaining-vf-cat}): as we shall see, tautness will play an essential role in providing this assignment with compositions. We will then show that this new vf-category is itself taut (\Cref{prop:FX-is-taut}), which will yield a proof of the desired result (\Cref{thm:main-theorem}).

\begin{rmk}
In this section, for convenience, we will always assume that a $\sigma$-family $(q_s)_{s:\sigma}$ of points in a vf-category $\X$ can be represented by some family $(q_s)_{s\in S}$ defined on the whole carrier set $S$ of the filter $\sigma$. Formally, this is not possible if (and only if) $\X$ is empty and $\sigma$ is the improper filter on a nonempty set $S$: in this case, we can replace $(S,\sigma)$ by the canonically isomorphic filter $0$ defined on the empty set. 
\end{rmk}

\subsection{Walking vf-arrows}\label{ssec:walking-vf-arrows}

In this subsection, we describe, for each filter $\sigma$ on a set $S$, a topological space $\mathrm I_{S,\sigma}$ which represents the formal convergence of $\sigma$-families of points. More concretely, the space $\mathrm{I}_{S,\sigma}$ is characterized by the property that vf-functors out of it correspond, in a suitable way, to vf-arrows of type $\sigma$ of the codomain vf-category. In the spirit of referring to the free categorical entity admitting some gadget $x$ as the \emph{walking $x$} --- e.g., the {walking arrow} $I \coloneqq \set{0 \to 1}$, which is such that functors $I \to \mathcal C$ correspond to an arrow in $\mathcal C$ --- we name these spaces \emph{walking vf-arrows}. 

Throughout this subsection, fix a filter $\sigma$ on a set $S$.

\begin{cons}[Walking vf-arrows]\label{cons:walking-vf-arrow}
\looseness=-1
    The \emph{walking vf-arrow} for $(S,\sigma)$ is the topological space $\mathrm I_{S,\sigma}$ whose set of points is given by $\set{ a } + \set{ b_s | s \in S}$ and whose opens are subsets $U$ such that 
    \[ a \in U \Rightarrow \forall_{s:\sigma} \ b_s \in U.\]
\end{cons}

\begin{rmk}
As our choice of notation emphasizes, the space $\mathrm{I}_{S,\sigma}$ explicitly depends on the carrier set $S$ of the filter $\sigma$; in particular, the assignment $(S,\sigma) \mapsto \mathrm{I}_{S,\sigma}$ does not extend to a functor $\FF \to \TopSp$. Although this might appear misaligned with our conventions in this paper, we shall see below that the key property characterizing $\mathrm{I}_{S,\sigma}$ depends only on a $\sigma$-large subset of $S$ (\Cref{cor:I-walking}).
\end{rmk}

By definition, seeing $\mathrm I_{S,\sigma}$ as a vf-category, we have that $a \vuleq{s:\sigma} b_s$ holds: we denote by $u_{S,\sigma} \colon a \vuto{s:\sigma} b_s$ the unique vf-arrow witnessing convergence. We can think of $\mathrm{I}_{S,\sigma}$ as the vf-category \emph{freely generated} by a vf-arrow of type $\sigma$: to substantiate this idea, we shall now describe the convergence structure of the space $\mathrm I_{S,\sigma}$ explicitly.

\begin{lem}[Convergence in walking vf-arrows]\label{lem:convergence-in-walking-vf-arrows}
    In $\mathrm I_{S,\sigma}$, vf-inequalities can be uniquely written in one of the following forms:
    \begin{enumerate}
        \item diagonals
        \[\diag^{b_s}_\lambda \colon b_s\vuleq{\lambda}b_s\]
        for some $s\in S$ and some filter $\lambda$;
        \item merges
        \[\diag^a_\kappa\mrg \phi^*u_{S,\sigma} \colon a\vuleq{\kappa+\tau}((a)_{\kappa},(b_{\phi t})_{\tau})\]
        for some $\FF$-arrow $\phi \colon \tau\to\sigma$ and some filter $\kappa$  (cf.\ \Cref{nota:merging-pairs}).
    \end{enumerate}
\end{lem}
\begin{proof}
    Consider a vf-inequality $x\vuleq{\lambda}c_\ell$ in $\mathrm I_{S,\sigma}$. If $x = b_s$ for some $s\in S$, then necessarily $c_\ell = b_s$ for $\ell:\lambda$, and we are in the first case. Now assume that $x = a$. Seeing the $\lambda$-family $(c_\ell)_{\ell:\lambda}$ as an $\FF$-arrow
\[\begin{tikzcd}
	{\lambda } & {\set{a} + \set{b_s | s\in S}}
	\arrow["c", from=1-1, to=1-2]
\end{tikzcd}\]
    note that, by extensivity of $\FF$ (\Cref{prop:F-is-lextensive}), we can decompose $c$ as $c_1 + c_2$ for appropriate $\FF$-arrows $c_1 \colon \kappa \to \set{ a }$ and $c_2 \colon \tau \to \set{ b_s | s \in S}$. By the description of the topology on $\mathrm I_{S,\sigma}$, from $a \vuleq{\lambda} c_\ell$ it follows that $c_2$ factors through $\sigma$; in other words, there exists an $\FF$-arrow $\phi \colon \tau \to \sigma$ such that $c_2$ factors as the composite
\[\begin{tikzcd}
	\tau & \sigma & {\set{b_s | s\in S}}
	\arrow["\phi", from=1-1, to=1-2]
	\arrow[hook, from=1-2, to=1-3]
\end{tikzcd}.\]
The $\lambda$-family $(c_\ell)_{\ell:\lambda}$ can thus be written as $((a)_{\kappa},(b_{\phi t})_{t:\tau})$, and we are in the second case. In either case, uniqueness follows from posetality.
\end{proof}

The previous lemma entails that a vf-functor $F$ out of the space $\mathrm{I}_{S,\sigma}$ is uniquely determined by its action on the vf-arrow $u_{S,\sigma}$, once the images $F(a)$ and $(F(b_s))_{s\in S}$ have been fixed. Before stating the `walking' property of $\mathrm{I}_{S,\sigma}$, intuitively asserting that the correspondence $F \mapsto F(u_{S,\sigma})$ is bijective, we need a technical observation: for a vf-category $\X$, the set of vf-functors $F \colon \mathrm{I}_{S,\sigma} \to \X$ mapping $a$ and each $b_s$ to prescribed points of $\X$ only depends, in fact, on a $\sigma$-large subset of the $b_s$'s. Let us introduce some notation. 

\begin{cons}
    Fixing a point $p$ and an $S$-family of points $(q_s)_{s\in S}$ allows us to see any vf-category $\X$ as a \emph{$(1+S)$-pointed} vf-category, that is, an object of the coslice 2-category $(1+S)\downarrow\vfCAT$. Seeing the space $\mathrm{I}_{S,\sigma}$ as $(1+S)$-pointed in the evident way, we can then consider the hom-category 
    \[ \X_\sigma [p, q_s]\coloneqq ((1+S)\downarrow \vfCAT)(\mathrm I_{S,\sigma},\X),\]
    whose objects are vf-functors $\mathrm{I}_{S,\sigma} \to \X$ mapping $a \mapsto p$ and $b_s \mapsto q_s$ for each $s\in S$. 

\end{cons}

The following proposition shows that we can consider $\X_{\sigma}[p,q_s]$ more generally for any point $p$ and any $\sigma$-family of points $(q_s)_{s:\sigma}$, that is, $\X_{\sigma}[p,q_s]$ only depends on a $\sigma$-large subset of the $(q_s)_{s\in S}$.

\begin{prop}\label{prop:cocomma-depends-on-large-set}
    Let $\X$ be a vf-category and fix a point $p$ and an $S$-family of points $(q_s)_{s\in S}$. Then, the category  $\X_\sigma [p, q_s] $ is small and discrete, i.e.\ a set. 
    Moreover, if $(q_s')_{s\in S}$ is another $S$-family of points such that $q_s = q_s'$ for $s :  \sigma$, then there is a canonical bijection:
    \[ \X_\sigma [p, q_s] \cong \X_\sigma [p, q_s'].\]
\end{prop}
\begin{proof}
That the category $\X_\sigma [p, q_s]$ is discrete is evident as the function $1+S \to \mathrm{I}_{S,\sigma}$ is the identity on points; moreover, it is small by \Cref{lem:convergence-in-walking-vf-arrows} since so is $\X_\sigma (p, q_s)$.
For the second part, note first that $\mathrm{I}_{S,\sigma}$ is the coproduct of $\mathrm{I}_{S_0,\sigma}$ and the discrete space $S \setminus S_0$, in $\TopSp$ and hence in $\vfCAT$\footnote{Indeed, it is straightforward to see that the embeddings $\TopSp \hookrightarrow \vfCAT$ and $\TopSp \hookrightarrow \vuCAT$ preserve small coproducts.}, where $S_0$ is any $\sigma$-large subset of $S$ and where we identify $\sigma$ with its restriction to $S_0$. This yields a canonical bijection
\begin{multline*}
    ((1+S)\downarrow \vfCAT)(\mathrm I_{S,\sigma},\X) \cong \\ ((1+S_0)\downarrow \vfCAT)(\mathrm I_{S_0,\sigma},\X) \times ((S\setminus S_0)\downarrow \vfCAT)(S\setminus S_0,\X) ,
\end{multline*}
where of course $((S\setminus S_0)\downarrow \vfCAT)(S\setminus S_0,\X)$ is a singleton. In other words, $\X_\sigma [p, q_s]$ only depends on a $\sigma$-large subset of the $q_s$'s, whence the claim.
\end{proof}

\begin{cor}[The `walking' property]\label{cor:I-walking}
    Let $\X$ be a vf-category and fix a point $p$ and a $\sigma$-family of points $(q_s)_{s: \sigma}$. Then, the 
    map
    \[ \X_\sigma[p, q_s] \to \X_\sigma (p, q_s) \qquad F \mapsto F(u_{S,\sigma})\]
    is a bijection.
\begin{proof}

    \looseness=-1
By \Cref{prop:cocomma-depends-on-large-set} we can fix any $S$-family $(q_s)_{s\in S}$ representing $(q_s)_{s: \sigma}$, so that elements of $\X_\sigma[p, q_s]$ are identified with vf-functors $F \colon \mathrm{I}_{S,\sigma} \to \X$ mapping $a \mapsto p$ and $b_s \mapsto q_s$ for each $s\in S$. Then, by \Cref{lem:convergence-in-walking-vf-arrows}, we clearly have that the map $F \mapsto F(u_{S,\sigma})$ is injective, since such a vf-functor is uniquely determined by its action on $u_{S,\sigma}$.

To show that the map is surjective, let $w \colon p \vuto{\sigma} q_s$ be a vf-arrow in $\X$. By \Cref{lem:convergence-in-walking-vf-arrows}, the assignment $a \mapsto p$, $b_s \mapsto q_s$ for each $s\in S$, and $u_{S,\sigma} \mapsto w$ uniquely extends to a vf-functor $F_w \colon \mathrm{I}_{S,\sigma} \to \X$. More explicitly, $F_w$ acts on vf-inequalities as follows:
\[ F_w(\diag^{b_s}_\kappa) = \diag^{q_s}_\kappa \qquad \text{and} \qquad F_w(\diag^a_\kappa \mrg \phi^* u_{S,\sigma}) = \diag^p_\kappa \mrg \phi^* w.\]
We omit the routine check that $F_w$ so defined is indeed a vf-functor.
\end{proof}
\end{cor}

\begin{rmk}[Walking vu-arrows]\label{rmk:I-walking-vu}
    The proofs of this subsection adapt straightforwardly to vu-categories. More precisely, for a vu-category $\X$, a point $p$, and a $\sigma$-family of points $(q_s)_{s:\sigma}$, the set 
    \[\X_\sigma [p, q_s] \coloneqq ((1+S)\downarrow \vuCAT)(\mathrm{I}_{S,\sigma}, \X)\] 
    does not depend on the $S$-family representing $(q_s)_{s:\sigma}$; moreover, if $\sigma$ is in fact an ultrafilter, the map $\X_\sigma [p, q_s] \to \X_\sigma (p, q_s)$ defined by $F\mapsto F(u_{S,\sigma})$ is a bijection, and thus $\mathrm{I}_{S,\sigma}$ acts as a \emph{walking vu-arrow}.
\end{rmk}

\subsection{From taut vu-categories to taut vf-categories}\label{ssec:from-vu-to-vf}Exploiting walking vf-arrows, we will construct a vf-category $\mathsf F \X$ out of a fixed vu-category $\X$, crucially assumed to be taut. The key idea is to define $\mathsf F \X$ as the vf-category whose points are those of $\X$ and whose set of vf-arrows $p \vuto{\sigma} q_s$, for a filter $\sigma$ on a set $S$, is given by the set $\X_\sigma [p, q_s]$ constructed in the previous subsection. Concretely, we will proceed as follows:
\begin{enumerate}
    \item[(a)]first, we show how to endow this assignment with a reindexing operation;
    \item[(b)]next, we show how to endow it with a composition operation;
    \item[(c)]then, we assemble this structure into a vf-category $\mathsf F \X$;
    \item[(d)]finally, we prove that $\mathsf F \X$ is taut.
\end{enumerate}

\subsubsection{Reindexings}\label{ssec:reindexing}

First, we deal with reindexings: fix an $\FF$-arrow $\phi \colon \tau \to \sigma$ represented by a total function $\phi \colon T \to S$.

\begin{nota}\label{nota:zeta-phi}
In $\mathrm I_{S,\sigma}$, the vf-inequality $a \vuleq{\tau} b_{\phi t}$ holds; thus, by \Cref{cor:I-walking}, it uniquely determines a continuous map 
\[ \zeta_\phi \colon \mathrm I_{T,\tau} \to \mathrm I_{S,\sigma}\]
acting as $a \mapsto a$ and $b_t \mapsto b_{\phi t}$ for each $t \in T$.  
\end{nota}

\begin{cons}[Reindexings]\label{cons:vu-to-vf-reindexings}
Fix an element $F \in \X_\sigma[p, q_s]$. Up to fixing some representatives $(q_s)_{s\in S}$, we can identify $F$ with a vu-functor $\mathrm I_{S,\sigma} \to \X$. Precomposing with $\zeta_\phi$ we obtain a vu-functor $\mathrm I_{T,\tau} \to \X$, representing an element of $\X_\tau[p, q_{\phi t}]$ which we denote by $\phi^* F$.
\end{cons}

\begin{lem}\label{lem:vu-to-vf-reind}
With the above notation, the map 
\[ \X_\sigma[p, q_s] \to \X_\tau[p, q_{\phi t}] \qquad F \mapsto \phi^*F\]
is well-defined, and it induces a functor
\[(\FF/\sigma)^{\op}\to\Set \qquad (\phi\colon \tau\to\sigma) \mapsto \X_\tau[p, q_{\phi t}].\]
\begin{proof}
    By \Cref{prop:cocomma-depends-on-large-set}. 
\end{proof}
\end{lem}

\subsubsection{Compositions}\label{ssec:composition}
In order to define compositions, we will `glue' walking vf-arrows appropriately: fix a filter $\sigma$ on a set $S$ and, for each $s \in S$, a filter $\tau_s$ on a set $T_s$. Denoting by $\set{ b_s} + \set{ c_{s,t} | t \in T_s}$ the set of points of $\mathrm{I}_{T_s, \tau_s}$, we will consider the space  
obtained by attaching $\mathrm{I}_{S,\sigma}$ to the coproduct $S \cdot \mathrm{I}_{T_s,\tau_s} \coloneqq \coprod_{s\in S} \mathrm{I}_{T_s,\tau_s}$ along the map $S\cdot b_s \colon S \to S \cdot \mathrm{I}_{T_s,\tau_s}$, where $S$ is seen as a discrete subspace of $\mathrm{I}_{S,\sigma}$. 

\begin{cons}[Walking vf-composition]
   With the above notation, we denote by $\mathrm J_{S,\sigma : T_s,\tau_s}$ the following pushout in $\TopSp$:
    \[\begin{tikzcd}[sep = small]
        S && {S\cdot \mathrm{I}_{T_s,\tau_s}} \\
        \\
        {\mathrm{I}_{S,\sigma}} &&  {\mathrm J_{S,\sigma : T_s,\tau_s}}
        \arrow["S\cdot b_s", from=1-1, to=1-3]
        \arrow["b"', hook, from=1-1, to=3-1]
        \arrow[from=1-3, to=3-3]
        \arrow[from=3-1, to=3-3]
        \arrow["\lrcorner"{anchor=center, pos=0.125, rotate=180}, draw=none, from=3-3, to=1-1]
    \end{tikzcd}\]
    Explicitly, $\mathrm J_{S,\sigma : T_s,\tau_s}$ is the space whose points are given by
    \[ \set{ a } + \set{ b_s | s \in S } + \set{ c_{s,t} | s \in S, t \in T_s },\]
    and whose opens are subsets $U$ such that 
    \[ a \in U \Rightarrow \forall_{s:\sigma} \ b_s\in U \quad \mbox{ and }\quad \forall s \in S \left(\ b_s \in U \Rightarrow \forall_{t:\tau_s} \ c_{s,t} \in U \  \right). \]
\end{cons}

\begin{lem}\label{lem:composition-pushout-preserved}
    The above pushout 
    is preserved by the inclusion $\TopSp\mono\vuCAT_{\mathrm{taut}}$.
\end{lem}
\begin{proof}
    The proof rests entirely on the fact that the pushout of $\mathrm{I}_{S,\sigma} \hookleftarrow S \to S \cdot \mathrm{I}_{T_s,\tau_s}$ computed in $\vuCAT$ exists and is, in fact, a vu-poset: to see this, note that we can describe this pushout $P$ explicitly as the vu-poset on the same points as $\mathrm J_{S,\sigma : T_s,\tau_s}$ whose vu-inequalities are freely generated by the vu-inequalities of $\mathrm{I}_{S,\sigma}$ and of $S \cdot \mathrm{I}_{T_s,\tau_s}$. 
 
    On the one hand, seeing $\mathrm J_{S,\sigma : T_s,\tau_s}$ as a taut vu-poset, the universal property of $P$ induces a vu-functor $F \colon P \to \mathrm J_{S,\sigma : T_s,\tau_s}$
    which, by the universal property of the tautification maps (\Cref{defn:tautification}), descends to a vu-functor $\widetilde F \colon \Tref{P} \to \mathrm J_{S,\sigma : T_s,\tau_s}$. On the other hand, as $\Tref{P}$ is a topological space, the universal property of $\mathrm J_{S,\sigma : T_s,\tau_s}$ induces a continuous map (i.e.\ a vu-functor) $G \colon \mathrm J_{S,\sigma : T_s,\tau_s} \to \Tref{P}$. Composing the two universal properties, it follows that $\widetilde F$ and $G$ are each other's inverse, so that $\Tref{P}$ and $\mathrm J_{S,\sigma : T_s,\tau_s}$ are homeomorphic (i.e.\ isomorphic as vu-posets). 

    By \Cref{prop:tautness-orthogonality}, we then have an isomorphism of categories
    \[ \vuCAT(\mathrm J_{S,\sigma : T_s,\tau_s}, \Y) \cong \vuCAT(P, \Y)\] 
    for any taut vu-category $\Y$; thus, the universal property of $P$ in $\vuCAT$ implies the analogous universal property for $\mathrm J_{S,\sigma : T_s,\tau_s}$ in $\vuCAT_{\mathrm{taut}}$. 
\end{proof}

\begin{rmk}[On the role of tautness]
    We do not know whether, in general, $\mathrm J_{S,\sigma : T_s,\tau_s}$ is the pushout of $\mathrm{I}_{S,\sigma} \hookleftarrow S \to S \cdot \mathrm{I}_{T_s,\tau_s}$ in $\vuCAT$, or equivalently whether the vu-poset $P$ constructed in the proof above is itself taut. Although we expect that there exist filters $\sigma$ and $(\tau_s)_{s\in S}$ for which this is not the case, we have not been able to construct any. It even seems plausible that the existence of such filters is independent of ZFC; in this picture, tautness appears as a shield against set-theoretic pathologies.
\end{rmk}

\begin{nota}
In $\mathrm J_{S,\sigma : T_s,\tau_s}$, the vf-inequality $a \vuleq{\sigma\cdot \tau_s} c_{s,t}$ holds; thus, by \Cref{cor:I-walking}, it uniquely determines a continuous map 
\[ \Theta_{\sigma\cdot \tau_s} \colon \mathrm I_{S\cdot T_s , \sigma \cdot \tau_s} \to \mathrm J_{S,\sigma : T_s,\tau_s}.\]
\end{nota}

\begin{cons}[Compositions]\label{cons:vu-to-vf-compositions}
    Fix an element $F \in \X_\sigma[p, q_s]$ and a $\sigma$-family of elements $(G_s \in \X_{\tau_s}[q_s, r_{s,t}])_{s:\sigma}$. Up to fixing some representatives $(q_s)_{s\in S}$ and $((r_{s,t})_{t\in T_s})_{s\in S}$, we can identify $F$ with a vu-functor $\mathrm I_{S,\sigma} \to \X$ and each $G_s$ with a  vu-functor $ \mathrm I_{T_s,\tau_s} \to \X$. 
\[\begin{tikzcd}
	& {S\cdot\mathrm I_{T_s,\tau_s}} & &&[14pt] \\
	{S } & {\mathrm I_{S\cdot T_s, \sigma\cdot\tau_s}} & {\mathrm J_{S,\sigma:T_s,\tau_s}} & & \X \\
	& {\mathrm I_{S,\sigma}}
	\arrow[from=1-2, to=2-3]
	\arrow["{(G_s)_{s\in S}}", curve={height=-24pt}, from=1-2, to=2-5]
	\arrow["{S\cdot b_s}", from=2-1, to=1-2]
	\arrow["{b_s}"', hook, from=2-1, to=3-2]
	\arrow["\Theta_{\sigma\cdot \tau_s}", from=2-2, to=2-3]
	\arrow["{((G_s)_{s\in S}, F)}", dashed, from=2-3, to=2-5]
	\arrow[from=3-2, to=2-3]
	\arrow["F"', curve={height=24pt}, from=3-2, to=2-5]
\end{tikzcd}\]
By the universal properties of coproducts and pushouts, we then have a canonical vu-functor $((G_s)_{s\in S}, F) \colon \mathrm J_{S,\sigma : T_s,\tau_s} \to \X$. Precomposing with $\Theta_{\sigma\cdot \tau_s}$ we obtain a vu-functor $\mathrm I_{S\cdot T_s, \sigma\cdot\tau_s} \to \X$, representing an element of $\X_{\sigma\cdot\tau_s}[p, r_{s,t}]$ which we denote by $(G_s)_{s:\sigma}\circ_\sigma F$.
\end{cons}

\begin{lem}\label{lem:vu-to-vf-compo}
    With the above notation, the map 
    \[\X_\sigma[p,q_s]\times\textstyle\int_{s:\sigma} \X_{\tau_s}[q_s,r_{s,t}]\to \X_{\sigma\cdot\tau_s}[p,r_{s,t}] \qquad (F,(G_s)_{s:\sigma}) \mapsto (G_s)_{s:\sigma} \circ_\sigma F \]
    is well-defined.
\begin{proof}
    By \Cref{prop:cocomma-depends-on-large-set}.
\end{proof}
\end{lem}

\subsubsection{Assembling}\label{ssec:assembling}\looseness=-1
We can now assemble reindexings and compositions into the structure of a vf-category.

\begin{Thm}[From vu-categories to vf-categories]\label{thm:obtaining-vf-cat}
    Let $\X$ be a taut vu-category. he sets $\X_{\sigma}[p, q_s]$ endow the class of points of $\X$ with the structure of a vf-category, denoted by $\mathsf F \X$, where:
    \begin{enumerate}
        \item identity vf-arrows are represented by constant vu-functors $\mathrm I_{1,1} \to \X$;
        \item reindexings are given by \Cref{cons:vu-to-vf-reindexings};
        \item compositions are given by \Cref{cons:vu-to-vf-compositions}.
    \end{enumerate}
\end{Thm}
\begin{proof}
    Concretely, we need to check the axioms of a vf-category (cf.\ \Cref{defn:vf-cats}); functoriality of reindexing is ensured by \Cref{lem:vu-to-vf-reind}, so we only have to show unitality and associativity of composition and compatibility between reindexings and compositions. For the sake of readability, exploiting \Cref{prop:cocomma-depends-on-large-set} and the well-definedness of reindexings and compositions, we will implicitly assume fixed representatives of all families of points, $\FF$-arrows, and vf-arrows involved in the following computations.

    \looseness=-1
    First, we show right unitality; left unitality is proved similarly. Let $F \in \X_\sigma[p, q_s]$ and consider the identity $\id_p \in \X_1[p,p]$, represented by the constant vu-functor $\mathrm{I}_{1,1} \to \X$ of value $p$. The equality $F \circ_{1} \id_p = F$ follows since, through the identification $1\cdot \mathrm{I}_{S,\sigma} \cong \mathrm{I}_{1\cdot S, 1\cdot \sigma}$, the composite $((F)_{*\in1}, \id_p) \circ \Theta_{1\cdot \sigma}$ coincides with $F$, as witnessed by the diagram below (where we denote the points of $\mathrm I_{1,1}$ by $\set{ a} + \set{b}$ and the points of $\mathrm I_{S,\sigma}$ by $\set{ b } + \set{ c_s | s \in S}$).
\[\begin{tikzcd}
	& {1\cdot\mathrm I_{S,\sigma}} &&&[14pt] \\
	1 & {\mathrm I_{1\cdot S, 1\cdot\sigma}} & {\mathrm J_{1,1:S,\sigma}} && \X \\
	& {\mathrm I_{1,1}}
	\arrow["\cong"', from=1-2, to=2-2]
	\arrow[from=1-2, to=2-3]
	\arrow["{(F)_{*\in 1}}", curve={height=-24pt}, from=1-2, to=2-5]
	\arrow["{1\cdot b}", from=2-1, to=1-2]
	\arrow["b"', hook, from=2-1, to=3-2]
	\arrow["{\Theta_{1\cdot \sigma}}", from=2-2, to=2-3]
	\arrow["{((F)_{*\in 1}, \id_p)}", dashed, from=2-3, to=2-5]
	\arrow[from=3-2, to=2-3]
	\arrow["{\id_p}"', curve={height=24pt}, from=3-2, to=2-5]
\end{tikzcd}\]

    Next, we show right compatibility of composition with reindexing; left compatibility is proved similarly. Let $F \in \X_\sigma[p, q_s]$, $(G_s)_{s:\sigma} \in \int_{s:\sigma} \X_{\tau_s}[q_s, r_{s,t}]$ and $\phi \colon \lambda \to \sigma$. Unraveling definitions, we need to show that the following diagram of vu-functors commutes:
\[(\star) \qquad \begin{tikzcd}[column sep = 30pt, row sep =20pt]
	{\mathrm I_{L\cdot T_{\phi \ell}, \lambda \cdot \tau_{\phi \ell}}} & {\mathrm I_{S\cdot T_s, \sigma\cdot\tau_s}} & {\mathrm J_{S,\sigma:T_s,\tau_s}} \\
	\\
	{\mathrm J_{L,\lambda:T_{\phi \ell},\tau_{\phi \ell}}} && \X
	\arrow["{\zeta_{\phi\cdot \operatorname{id}_{\tau_s}}}", from=1-1, to=1-2]
	\arrow["{\Theta_{\lambda\cdot \tau_{\phi \ell}}}"', from=1-1, to=3-1]
	\arrow["{\Theta_{\sigma\cdot \tau_s}}", from=1-2, to=1-3]
	\arrow["{((G_s)_{s\in S}, F)}", dashed, from=1-3, to=3-3]
	\arrow["{((G_{\phi \ell})_{\ell\in L}, F\zeta_\phi )}"', dashed, from=3-1, to=3-3]
\end{tikzcd}\]
Note then that the diagram below, where we denote by $\set{ a' } + \set{ b'_\ell | \ell \in L}$ the points of $\mathrm I_{L,\lambda}$, commutes. Thus, by the universal property of the pushout, we uniquely obtain a vu-functor $\Omega \colon \mathrm J_{L,\lambda:T_{\phi \ell},\tau_{\phi \ell}} \to \mathrm J_{S,\sigma:T_s,\tau_s}$ commuting with the two dashed arrows in $(\star)$.
\[\begin{tikzcd}[column sep = 30pt, row sep =15pt]
	{\mathrm I_{L,\lambda}} &[10pt] L & {L\cdot\mathrm I_{T_{\phi \ell},\tau_{\phi \ell}}} \\
	\\
	{\mathrm I_{S,\sigma}} & {S } & {S\cdot\mathrm I_{T_s,\tau_s}}
	\arrow["{\zeta_\phi}"', from=1-1, to=3-1]
	\arrow["{b'_\ell}"', hook', from=1-2, to=1-1]
	\arrow["{L\cdot b_{\phi \ell}}", from=1-2, to=1-3]
	\arrow["\phi"', from=1-2, to=3-2]
	\arrow["{(\iota_{\phi \ell})_{\ell \in L}}", from=1-3, to=3-3]
	\arrow["{b_s}", hook', from=3-2, to=3-1]
	\arrow["{S\cdot b_s}"', from=3-2, to=3-3]
\end{tikzcd}\]
Explicitly, $\Omega$ acts on the points of $\mathrm J_{L,\lambda:T_{\phi \ell},\tau_{\phi \ell}}$, which we denote by $\set{ a' } + \set{ b'_\ell | \ell \in L} + \set{ c'_{\ell, t} | \ell\in L, t \in T_{\phi \ell}}$, as follows:
\[ a' \mapsto a \qquad b'_\ell \mapsto b_{\phi \ell} \qquad c'_{\ell, t} \mapsto c_{\phi \ell, t}.\]
Thus, by direct inspection and since a vu-functor between vu-posets is uniquely determined by its action on the points, $\Omega$ also commutes with the two solid arrows in $(\star)$; hence, the whole diagram $(\star)$ commutes.

Associativity is proved in a similar --- albeit more tedious --- way, which we omit.
\end{proof}

\subsubsection{Tautness}\label{ssec:FX-is-taut}\looseness=-1 Finally, we now show that the vf-category $\mathsf F \X$ is taut: to do this, we need another technical result similar to the above \Cref{lem:composition-pushout-preserved}. Throughout, let $\phi \colon \tau \epi \sigma$ be a cover in $\FF$ with a fixed representative, given by a total surjection $\phi \colon T \epi S$.

\begin{cons}\label{cons:walking-phi-reduc}
    We denote by $\mathrm I_{S,\phi}$ the vf-poset on the set $\set{ a } + \set{ b_s | s \in S }$ generated by the vf-inequality $a \vuleq{\tau} b_{\phi t}$: that is, the smallest vf-poset obtained from $a \vuleq{\tau} b_{\phi t}$ by reindexings and compositions.
\end{cons}

\begin{rmk}\label{rmk:I-S-phi-sub-vf-poset}
    The identity map on points extends to a vf-functor $\mathrm I_{S,\phi} \to \mathrm I_{S,\sigma}$, so that we can identify $\mathrm I_{S,\phi}$ with a \emph{sub-vf-poset} of $\mathrm I_{S,\sigma}$ on the same points. Through this identification, in the spirit of \Cref{lem:convergence-in-walking-vf-arrows}, note that vf-inequalities in $\mathrm I_{S,\phi}$ are either diagonals at some $b_s$ or of the form $a\vuleq{\kappa+\lambda}((a)_\kappa,(b_{\psi \ell}))$ for some $\FF$-arrow $\psi\colon \lambda\to\sigma$ factoring through $\phi\colon \tau\epi\sigma$. 
\end{rmk}

\looseness=-1 Recall from \Cref{nota:zeta-phi} that the vf-functor $\zeta_\phi \colon \mathrm I_{T,\tau} \to \mathrm I_{S,\sigma}$ corresponds to the vf-arrow $\phi^*u_{S,\sigma} \colon a \vuleq{\tau} b_{\phi t}$ in $\mathrm I_{S,\sigma}$; as $\phi^*u_{S,\sigma}$ lies in $\mathrm I_{S,\phi}\subseteq \mathrm I_{S,\sigma}$, the vf-functor $\zeta_{\phi}$  factors through $\mathrm I_{S,\phi}$. We then have the following property of $\mathrm I_{S,\phi}$.

\begin{lem}
    The diagram below is a coequalizer in $\vfCAT$.
    \[\begin{tikzcd}
        {\mathrm{I}_{T \times_S T , \tau \times_\sigma \tau}} & {\mathrm I_{T,\tau}} & {\mathrm I_{S,\phi}}
        \arrow["{\zeta_{\pi_1}}", shift left, from=1-1, to=1-2]
        \arrow["{\zeta_{\pi_2}}"', shift right, from=1-1, to=1-2]
        \arrow["{\zeta_{\phi}}", from=1-2, to=1-3]
    \end{tikzcd}\] 
\end{lem}
\begin{proof}
    The proof is similar to that of the next lemma (\Cref{lem:I-S-phi-is-coequalizer}).
\end{proof}

\begin{rmk}[Walking reducible vf-arrows]\label{rmk:walking-phi-reduc}
   \looseness=-1 Note then that, by \Cref{cor:I-walking}, a vf-functor out of $\mathrm I_{T,\tau}$ coequalizing $\zeta_{\pi_1}$ and $\zeta_{\pi_2}$ corresponds to a \emph{$\phi$-reducible} vf-arrow in the target vf-category: therefore, we can think of $\mathrm I_{S,\phi}$ as the \emph{walking $\phi$-reducible vf-arrow}. 
\end{rmk}

In the rest of this section we will not use the above lemma, but rather a similar property of the vu-poset underlying $\mathrm I_{S,\phi}$, proved below.
Note that the vu-inequalities of $\mathsf U \mathrm I_{S,\phi}$ are either diagonals or of the form $a\vuleq{\lambda}b_{\psi \ell}$ for some $\FF$-arrow $\psi\colon\lambda\to\sigma$, with $\lambda$ an ultrafilter, factoring through $\phi\colon\tau\epi\sigma$; this follows from \Cref{rmk:I-S-phi-sub-vf-poset} noting that $\kappa + \lambda$ is an ultrafilter only if either $\kappa$ or $\lambda$ is $0$. 

\begin{lem}\label{lem:I-S-phi-is-coequalizer}
    The diagram below is a coequalizer in $\vuCAT$.
    \[\begin{tikzcd}
        {\mathrm{I}_{T \times_S T , \tau \times_\sigma \tau}} & {\mathrm I_{T,\tau}} & {\mathsf U\mathrm I_{S,\phi}}
        \arrow["{\zeta_{\pi_1}}", shift left, from=1-1, to=1-2]
        \arrow["{\zeta_{\pi_2}}"', shift right, from=1-1, to=1-2]
        \arrow["{\zeta_{\phi}}", from=1-2, to=1-3]
    \end{tikzcd}\]
\begin{proof}

Fix a vu-category $\Z$ and a vu-functor  $F \colon \mathrm I_{T,\tau} \to \Z$ coequalizing $\zeta_{\pi_1}$ and $\zeta_{\pi_2}$; we shall show that $F$ descends to $\widetilde F \colon \mathsf U\mathrm I_{S,\phi} \to \Z$. Clearly, the assignment on points uniquely descends to $\mathrm{I}_{S,\phi}$ since the underlying diagram of sets and functions is a coequalizer, which defines $\widetilde F$ on points. Now let $v\colon a\vuleq{\lambda}b_{\psi \ell}$ be a vu-inequality in $\mathsf U\mathrm I_{S,\phi}$, for some $\FF$-arrow $\psi\colon\lambda\to\sigma$ factoring through $\phi\colon\tau\epi\sigma$ with $\lambda$ an ultrafilter. Take an $\FF$-arrow $\chi\colon\lambda\to\tau$ such that $\psi = \phi \chi$ and we set
\[\widetilde F(v)\coloneqq F(\chi^*u_{T,\tau}).\]
Since $F\zeta_{\pi_1} = F\zeta_{\pi_2}$, note that this definition does not depend on the choice of $\chi$. Indeed, for another such factorization $\chi'$, we can consider $(\chi,\chi')\colon\lambda\to\tau\times_\sigma\tau$ and we have:
\begin{align*}
    F(\chi^*u_{T,\tau}) & = F((\chi,\chi')^*(\pi_1^*u_{T,\tau}))\\
    & = F((\chi,\chi')^*(\zeta_{\pi_1}(u_{T\times_S T,\tau\times_\sigma\tau})))\\
    & = F\zeta_{\pi_1}((\chi,\chi')^*u_{T\times_S T,\tau\times_\sigma\tau})
\end{align*}
and similarly for $F((\chi')^*u_{T,\tau})$, so we can deduce from $F\zeta_1 = F\zeta_2$ that $F(\chi^*u_{T,\tau}) = F((\chi')^*u_{T,\tau})$.
With this definition, it is then routine to check that $\widetilde F \colon \mathsf U\mathrm I_{S,\phi} \to \Z$ is a vu-functor and that it is the unique one extending $F$.
\end{proof}    
\end{lem}

The previous lemma yields the property we will need to show that the vf-category $\mathsf F \X$ is taut. Indeed, note that the following diagram is a coequalizer in $\TopSp$:
\[\begin{tikzcd}
        {\mathrm{I}_{T \times_S T , \tau \times_\sigma \tau}} & {\mathrm I_{T,\tau}} & {\mathrm I_{S,\sigma}}
        \arrow["{\zeta_{\pi_1}}", shift left, from=1-1, to=1-2]
        \arrow["{\zeta_{\pi_2}}"', shift right, from=1-1, to=1-2]
        \arrow["{\zeta_{\phi}}", from=1-2, to=1-3]
    \end{tikzcd}\]
Then, \Cref{lem:I-S-phi-is-coequalizer} guarantees that the coequalizer of $\zeta_{\pi_1}$ and $\zeta_{\pi_2}$ in $\vuCAT$ exists and is, in fact, a vu-poset: proceeding exactly as in the proof of \Cref{lem:composition-pushout-preserved}, we thus obtain the following.

\begin{cor}\label{cor:lem:I-S-sigma-is-coequalizer}
   The above coequalizer is preserved by the inclusion $\TopSp \hookrightarrow \vuCAT_{\mathrm{taut}}$.
   
\end{cor}

\begin{rmk}[Open surjections]
The map $\zeta_\phi\colon \mathrm I_{T,\tau}\to\mathrm I_{S,\sigma}$ is an \emph{open surjection} of topological spaces, hence an effective epimorphism in $\TopSp$; the above results then show that $\zeta_\phi$ is also an effective epimorphism in $\vuCAT_{\mathrm{taut}}$, with a crucial role played by tautness. In fact, we believe that any open surjection of topological spaces is an effective epimorphism in $\vuCAT_{\mathrm{taut}}$, but the proof goes beyond the scope of this paper. This sheds some light on \cite[Rem.\ 7.14]{saadiaExtendingConceptualCompleteness2026}, and it suggests a version of \cite[Prop.\ 7.13]{saadiaExtendingConceptualCompleteness2026} for taut vu-categories where the lifting is assumed to exist only up to thickening.

\end{rmk}

\begin{prop}\label{prop:FX-is-taut}
    For a taut vu-category $\X$, the vf-category $\mathsf F\X$ is taut.
\begin{proof}
    Suppose that $F \in \X_\tau [p, q_{\phi t}]$ is $\phi$-reducible for some cover $\phi \colon \tau \epi \sigma$, meaning that $F \zeta_{\pi_1} = F \zeta_{\pi_2}$ where $\pi_1, \pi_2 \colon \tau \times_\sigma \tau \rightrightarrows \tau$ is the kernel pair of $\phi$. Showing that $F$ can be uniquely reduced amounts to showing that, in the diagram below, there exists a unique vu-functor $\mathrm{I}_{S,\sigma} \to \X$ making the triangle commute.
    \[\begin{tikzcd}
    	{\mathrm I_{T \times_S T,\tau\times_\sigma \tau}} & {\mathrm I_{T, \tau}} & {\mathrm{I}_{S,\sigma}} \\
    	&& \X
    	\arrow["{\zeta_{\pi_1}}", shift left, from=1-1, to=1-2]
    	\arrow["{\zeta_{\pi_2}}"', shift right, from=1-1, to=1-2]
    	\arrow["{\zeta_\phi}", from=1-2, to=1-3]
    	\arrow["F"', from=1-2, to=2-3]
    	\arrow[dashed, from=1-3, to=2-3]
    \end{tikzcd}\]
    This follows since, by \Cref{cor:lem:I-S-sigma-is-coequalizer}, the map $\zeta_{\phi} \colon \mathrm I_{T,\tau} \to \mathrm I_{S,\sigma}$ is the coequalizer of $\zeta_{\pi_1}$ and $\zeta_{\pi_2}$ in $\vuCAT_{\mathrm{taut}}$. 
\end{proof}
\end{prop}

\subsection{Proving the equivalence}\label{ssec:main-theorem}
In this subsection, we conclude the proof of our main result. Concretely, we will show that $\mathsf{FU}\Y \cong \Y$ for any taut vf-category $\Y$, i.e.\ that a taut vf-category can be recovered from its underlying vu-category. We will then collect all the pieces together in the final result (\Cref{thm:main-theorem}).

\looseness=-1 Let $\Y$ be a vf-category and let $\sigma$ be a filter on a set $S$. By the constructions of the previous subsections, for a point $p$ and a $\sigma$-family of points $(q_s)_{s:\sigma}$ in $\Y$, we can consider both sets 
\[ \Y_\sigma [p, q_s] \qquad \text{and} \qquad (\mathsf U \Y)_\sigma[p, q_s],\]
and in particular we have a canonical map $\Y_\sigma [p, q_s] \to (\mathsf U \Y)_\sigma[p, q_s]$ sending a vf-functor $F \colon \mathrm{I}_{S,\sigma} \to \Y$ to the vu-functor $\mathsf U F \colon \mathrm{I}_{S,\sigma} \to \mathsf U \Y$. 
As we now show, in \Cref{prop:hom-UF=id}, this map is bijective if $\Y$ is taut: through the bijection $\Y_\sigma [p, q_s] \cong \Y_\sigma (p, q_s)$ of \Cref{cor:I-walking}, this implies that vf-arrows of type $\sigma$ in $\Y$ can be recovered from the vu-functors out of $\mathrm{I}_{S,\sigma}$ and into the underlying vu-category $\mathsf U \Y$. 

The key technical ingredient is the following lemma, crucially relying on the ultrafilter principle (\Cref{lem:cover-by-uf}).

\begin{lem}\label{lem:I-vu-vf}
    Let $\Y$ be a taut vf-category and, for a filter $\sigma$ on a set $S$, let $F \colon \mathrm I_{S,\sigma} \to \mathsf U \Y$ be a vu-functor. Then, there exists a unique vf-arrow $w \colon F(a) \vuto{\sigma} F(b_s)$ in $\Y$ such that ${\iota}^* w = F(\iota^* u_{S,\sigma})$ for all $\iota \colon \tau \mono \sigma$ with $\tau$ an ultrafilter.
\begin{proof}
    Let $L$ be a set indexing the ultrafilters $(\iota_\ell \colon \sigma_\ell \mono \sigma)_{\ell\in L}$ refining $\sigma$, so that by \Cref{lem:cover-by-uf} the canonical $\FF$-arrow $\phi \colon L\cdot\sigma_\ell\epi\sigma$ is a cover. For each $\ell \in L$, denote by $u_{S,\sigma_\ell}$ the vf-arrow $\iota_\ell^* u_{S,\sigma}$ witnessing that $a \vuleq{\sigma_{\ell}}b_s$. Consider then the merge
    \[\textstyle v \coloneqq \bigmrg_{\ell : L} F(u_{S,\sigma_\ell}) \colon F(a) \vuto{L\cdot\sigma_\ell} F(b_s).\]

    If we show that $v$ is reducible along $\phi \colon L\cdot\sigma_\ell\epi\sigma$, then by tautness of $\Y$ it will follow that there exists a unique vf-arrow $w \colon F(a) \vuto{\sigma} F(b_s)$ such that $\phi^* w = v$, which by construction satisfies 
    \[ \iota_\ell^* w = (\sigma_\ell \hookrightarrow L\cdot \sigma_\ell)^* \phi^* w = (\sigma_\ell \hookrightarrow L\cdot \sigma_\ell)^* v = F(\iota_\ell^* u_{S,\sigma})\]
    for each $\ell \in L$. Uniqueness of $w$ will then follow from \Cref{cor:vureduc-ultrafilters}. Now let $\psi_1,\psi_2 \colon \kappa\rightrightarrows L\cdot\sigma_\ell$ be two $\FF$-arrows such that $\phi \psi_1 = \phi \psi_2$; we need to show that $\psi_1^* v = \psi_2^* v$, and by \Cref{prop:vureduc-vfreduc} we can assume that $\kappa$ is an ultrafilter.
    
    Let $\lambda_1$ be the image of the composite $\kappa \stackrel{\psi_1} \to L \cdot \sigma_\ell \to L$, which is an ultrafilter on $L$ since so is $\kappa$. Note that the pullback of $\lambda_1 \mono L$ along the projection $L\cdot \sigma_\ell \to L$ is the ultrafilter $\lambda_1 \cdot \sigma_\ell$ refining $L \cdot \sigma_\ell$, so we can let $\theta_1 \colon \kappa \to \lambda_1 \cdot \sigma_\ell$ be the induced $\FF$-arrow.\vspace{-4pt}
    \[\begin{tikzcd}
    	\kappa & {\lambda_1\cdot\sigma_\ell} & {L\cdot\sigma_\ell} \\
    	& {\lambda_1} & L
    	\arrow["{\theta_1}", dashed, from=1-1, to=1-2]
    	\arrow["{\psi_1}", curve={height=-30pt}, from=1-1, to=1-3]
    	\arrow[two heads, from=1-1, to=2-2]
    	\arrow["{j_1}", hook, from=1-2, to=1-3]
    	\arrow[from=1-2, to=2-2]
    	\arrow["\lrcorner"{anchor=center, pos=0.125}, draw=none, from=1-2, to=2-3]
    	\arrow[from=1-3, to=2-3]
    	\arrow[hook, from=2-2, to=2-3]
    \end{tikzcd}\]
    Thus, $\psi_1$ factors as $\theta_1$ followed by the evident monomorphism $j_1 \colon \lambda_1 \cdot \sigma_\ell \mono L\cdot\sigma_\ell$. Similarly, we can let $\lambda_2$ be the image of the composite $\kappa \stackrel{\psi_2} \to L \cdot \sigma_\ell \to L$ and factor $\psi_2$ as a universally-induced $\theta_2 \colon \kappa \to \lambda_2 \cdot \sigma_\ell$ followed by $j_2 \colon \lambda_2 \cdot \sigma_\ell \mono L\cdot\sigma_\ell$, so that the two composites below coincide.\vspace{-4pt}
\[\begin{tikzcd}[sep = small]
	& {\lambda_1\cdot\sigma_\ell} && \\
	\kappa && {L\cdot\sigma_\ell} & \sigma \\
	& {\lambda_2\cdot\sigma_\ell}
	\arrow["{j_1}", hook, from=1-2, to=2-3]
	\arrow["{\theta_1}", from=2-1, to=1-2]
	\arrow["{\theta_2}"', from=2-1, to=3-2]
	\arrow["\phi", two heads, from=2-3, to=2-4]
	\arrow["{j_2}"', hook', from=3-2, to=2-3]
\end{tikzcd}\vspace{-8pt}\]
We can then compute:
\[\textstyle \psi_1^*v  = \theta_1^*\left(j_1^* \left( \bigmrg_{\ell:L} F(u_{S,\sigma_\ell})\right) \right)
    = \theta_1^* \left( \bigmrg_{\ell:\lambda_1} F(u_{S,\sigma_\ell}) \right)
     = F\left(\theta_1^*\left( \bigmrg_{\ell:\lambda_1} u_{S,\sigma_\ell}\right) \right)\]
where $\bigmrg_{\ell:\lambda_1} u_{S,\sigma_\ell}$ witnesses $a \vuleq{(\ell, s) : \lambda_1 \cdot \sigma_\ell} b_s$ in $\mathrm I_{S,\sigma}$, and similarly for $\psi_2^* v$. Since $\phi \psi_1 = \phi \psi_2$, we then have that the two vf-inequalities 
\[\textstyle\theta_1^*\left( \bigmrg_{\ell:\lambda_1} u_{S,\sigma_\ell} \right) \colon a\vuleq{\kappa}b_{\phi \psi_1 k} \quad \text{ and }\quad 
    \theta_2^*\left( \bigmrg_{\ell:\lambda_2} u_{S,\sigma_\ell} \right) \colon a\vuleq{\kappa}b_{\phi \psi_2 k}\] 
coincide as $\mathrm I_{S,\sigma}$ is posetal; hence, we conclude that $\psi_1^* v = \psi_2^* v$.
    \end{proof}
\end{lem}

\begin{prop}[Recovering vf-arrows from the underlying vu-category]\label{prop:hom-UF=id}
    Let $\Y$ be a taut vf-category and fix a point $p$ and a $\sigma$-family of points $(q_s)_{s:\sigma}$. Then, the map
    \[\Y_\sigma[p,q_s]\to(\mathsf U \Y)_\sigma[p,q_s] \qquad F \mapsto \mathsf UF\]
 is a bijection.
\begin{proof}
        Fix an $S$-family $(q_s)_{s\in S}$ representing $(q_s)_{s:\sigma}$ and take a vu-functor $F \colon \mathrm I_{S,\sigma} \to \mathsf U \Y$ mapping $a \mapsto p$ and $b_s \mapsto q_s$ for each $s\in S$. By \Cref{lem:I-vu-vf}, there exists a unique vf-arrow $w \colon p \vuto{\sigma} q_s$ in $\Y$ such that $\iota^* w = F(\iota^* u_{S,\sigma})$ for each $\iota \colon \tau \hookrightarrow \sigma$ with $\tau$ an ultrafilter. By \Cref{cor:I-walking}, $w$ corresponds to a vf-functor $F_w \colon \mathrm I_{S,\sigma} \to \Y$ mapping $u_{S,\sigma} \mapsto w$, and by construction it is the unique one such that $\mathsf UF_w = F$.
\end{proof}
\end{prop}

Since by definition $(\mathsf{FU}\Y)_\sigma (p, q_s) \coloneqq (\mathsf U \Y)_\sigma[p, q_s]$, the bijections established in the above proposition immediately determine an isomorphism of vf-categories $\mathsf{FU}\Y \cong \Y$, for any taut vf-category $\Y$. This is the last piece of the puzzle needed for the proof of our main theorem. Indeed, our strategy will be to realize the construction $\X \mapsto \mathsf F\X $ as a left 2-adjoint to $\mathsf U \colon \vfCAT_{\mathrm{taut}}\to\vuCAT_{\mathrm{taut}}$: 
    \[\begin{tikzcd}
    	& {} & \\
	{\vfCAT_{\mathrm{taut}}} & {} & {\vuCAT_{\mathrm{taut}}.}
	\arrow["\dashv"{anchor=center, rotate=-90, pos=0.6}, draw=none, from=1-2, to=2-2]
	\arrow["{\mathsf U}"', from=2-1, to=2-3]
	\arrow["{\mathsf F}"', curve={height=30pt}, dashed, from=2-3, to=2-1]
    \end{tikzcd}\]
At that point, the fact that $\operatorname{id} \cong \mathsf {UF}$ follows simply from the walking property of walking vu-arrows, while the fact that $\mathsf{FU} \cong \operatorname{id}$ is expressed by \Cref{prop:hom-UF=id} above. 

\begin{Thm}[The main equivalence]\label{thm:main-theorem}
    The 2-functor $\mathsf U \colon \vfCAT_{\mathrm{taut}}\to\vuCAT_{\mathrm{taut}}$ is a 2-equivalence of 2-categories.
\end{Thm}

\begin{proof}
    First, let $\X$ be a taut vu-category and consider $\mathsf {UF} \X$. By construction, a vu-arrow $p \vuto{\sigma} q_s$ in $\mathsf {UF} \X$ is represented by a vu-functor $\mathrm I_{S,\sigma} \to \X$, which by \Cref{rmk:I-walking-vu} uniquely corresponds to a vu-arrow $p \vuto{\sigma} q_s$ in $\X$. Thus, the identity on points extends to an evident vu-functor $\eta_{\X} \colon \X \to \mathsf{UF}\X$ which is an isomorphism of vu-categories.

\looseness=-1
Now let $G \colon \X \to \mathsf U \Y$ be a vu-functor for some taut vf-category $\Y$; we shall construct a vf-functor $\bar G \colon \mathsf F \X \to \Y$ satisfying $(\mathsf U \bar G) \circ \eta_{\X} = G$. For a filter $\sigma$ on a set $S$, fix some $F \in \X_\sigma[p, q_s]$ represented by a vu-functor $F \colon \mathrm I_{S,\sigma} \to \X$. Composing with $G$ we obtain a vu-functor $G F \colon \mathrm I_{S,\sigma} \to \mathsf U \Y$ which, by \Cref{lem:I-vu-vf}, uniquely determines a vf-arrow $\bar G (F) \colon G(p) \vuto{\sigma} G(q_s)$ in $\Y$: the assignment $F \mapsto \bar G (F) $ extends the map $p \mapsto G(p)$ on points to the desired vf-functor $\bar G \colon \mathsf F \X \to \Y$. Note also that the requirement $(\mathsf U \bar G) \circ \eta_{\X} = G$ determines $\bar G$ uniquely.

The correspondence $G \mapsto \bar G$ thus extends straightforwardly to a bijection on 2-cells, yielding a natural isomorphism of categories 
    \[ \vfCAT_{\mathrm{taut}} (\mathsf F\X, \Y)\cong  \vuCAT_{\mathrm{taut}} (\X, \mathsf U\Y) .\]
    In other words, the assignment $\X \mapsto \mathsf F \X$ extends to a 2-functor $\mathsf F \colon \vuCAT_{\mathrm{taut}} \to \vfCAT_{\mathrm{taut}}$ which is left 2-adjoint to $\mathsf U$. Note then that $\eta_{\X}$ realizes an isomorphism $\X \cong \mathsf {UF} \X$ for any taut vu-category $\X$, while $\mathsf{FU}\Y \cong \Y$ for any taut vf-category $\Y$ follows from \Cref{prop:hom-UF=id}.
\end{proof}

\begin{rmk}[Constructive reconstruction]
    \looseness=-1 As a final remark, we can immediately deduce from the main result of \cite{saadiaExtendingConceptualCompleteness2026, hamadGeneralisedUltracategoriesConceptual2025, vangoolToposesEnoughPoints2026} an analogous reconstruction theorem for toposes with enough points in terms of their vf-categories of points described in \Cref{ex:points-of-toposes}: if $\mathcal E$ is such a topos, then $\mathcal E \simeq \vfCAT(\pt(\mathcal E), \Set)$. As this proof relies on vu-categories, it is clearly \emph{not constructive}: however, we believe vf-categories to be the first step towards a fully constructive understanding of the reconstruction theorem. We leave the matter to future research.
\end{rmk}

\printbibliography

\end{document}